\documentclass[3p,preprint,12pt]{elsarticle}

\usepackage{titlecaps}
\Addlcwords{the of into via for and of on in an to hp-finite with}

\journal{Journal Name}

\usepackage{tabulary,xcolor}
\usepackage{amsfonts,amsmath,amssymb}
\usepackage[normalem]{ulem}

\usepackage[T1]{fontenc}

\makeatletter
\let\save@ps@pprintTitle\ps@pprintTitle
\def\ps@pprintTitle{\save@ps@pprintTitle\gdef\@oddfoot{\footnotesize\itshape \null\hfill\today}}
\def\hlinewd#1{%
  \noalign{\ifnum0=`}\fi\hrule \@height #1%
  \futurelet\reserved@a\@xhline}

\AtBeginDocument{\ifNAT@numbers \biboptions{sort&compress}\fi}
\makeatother

\usepackage{ifluatex}
\ifluatex
\usepackage{fontspec}
\defaultfontfeatures{Ligatures=TeX}
\usepackage[]{unicode-math}
\unimathsetup{math-style=TeX}
\else 
\usepackage[utf8]{inputenc}
\fi 
\ifluatex\else\usepackage{stmaryrd}\fi

\usepackage{url,multirow,morefloats,floatflt,cancel}

\usepackage{makecell}
\makeatletter
\AtBeginDocument{\@ifpackageloaded{textcomp}{}{\usepackage{textcomp}}}
\makeatother

\usepackage{colortbl}
\usepackage{xcolor}
\usepackage{pifont}
\usepackage[nointegrals]{wasysym}
\makeatletter
\def\mcWidth#1{\csname TY@F#1\endcsname+\tabcolsep}

\def\cAlignHack{\rightskip\@flushglue\leftskip\@flushglue\parindent\z@\parfillskip\z@skip}
\def\rAlignHack{\rightskip\z@skip\leftskip\@flushglue \parindent\z@\parfillskip\z@skip}

\if@twocolumn\usepackage{dblfloatfix}\fi 
\AtBeginDocument{
\expandafter\ifx\csname eqalign\endcsname\relax
\def\eqalign#1{\null\vcenter{\def\\{\cr}\openup\jot\m@th
  \ialign{\strut$\displaystyle{##}$\hfil&$\displaystyle{{}##}$\hfil
      \crcr#1\crcr}}\,}
\fi
}

\let\lt=<
\let\gt=>
\def\processVert{\ifmmode|\else\textbar\fi}

\@ifundefined{subparagraph}{
\def\subparagraph{\@startsection{paragraph}{5}{2\parindent}{0ex plus 0.1ex minus 0.1ex}%
{0ex}{\normalfont\small\itshape}}%
}{}

\newcommand\role[1]{\unskip}
\newcommand\aucollab[1]{\unskip}
  
\@ifundefined{tsGraphicsScaleX}{\gdef\tsGraphicsScaleX{1}}{}
\@ifundefined{tsGraphicsScaleY}{\gdef\tsGraphicsScaleY{.9}}{}
\def\checkGraphicsWidth{\ifdim\Gin@nat@width>\linewidth
    \tsGraphicsScaleX\linewidth\else\Gin@nat@width\fi}

\def\checkGraphicsHeight{\ifdim\Gin@nat@height>.9\textheight
    \tsGraphicsScaleY\textheight\else\Gin@nat@height\fi}

\def\fixFloatSize#1{}
\let\ts@includegraphics\includegraphics

\def\inlinegraphic[#1]#2{{\edef\@tempa{#1}\edef\baseline@shift{\ifx\@tempa\@empty0\else#1\fi}\edef\tempZ{\the\numexpr(\numexpr(\baseline@shift*\f@size/100))}\protect\raisebox{\tempZ pt}{\ts@includegraphics{#2}}}}

\AtBeginDocument{\def\includegraphics{\@ifnextchar[{\ts@includegraphics}{\ts@includegraphics[width=\checkGraphicsWidth,height=\checkGraphicsHeight,keepaspectratio]}}}

\def\URL#1#2{\@ifundefined{href}{#2}{\href{#1}{#2}}}

\def\UrlOrds{\do\*\do\-\do\~\do\'\do\"\do\-}%
\g@addto@macro{\UrlBreaks}{\UrlOrds}
\makeatother

\usepackage{floatrow}

\usepackage{pgfplots,algorithmic,algorithm}
\pgfplotsset{compat=newest}

\usepackage{booktabs}

\usepackage[frozencache]{minted}

\usepackage[toc,page]{appendix}
\usepackage{bm}

\paperwidth=\dimexpr\paperwidth + 6cm\relax
\oddsidemargin=\dimexpr\oddsidemargin + 3cm\relax
\evensidemargin=\dimexpr\evensidemargin + 3cm\relax
\marginparwidth=\dimexpr\marginparwidth + 3cm\relax

\usepackage[textsize=scriptsize]{todonotes}
\newcommand\ed[1]{} 

\usepackage{pifont}
\newcommand{\cmark}{\ding{51}}
\newcommand{\xmark}{\ding{55}}

\usepackage{graphicx,amssymb,amsmath,amsthm}

\newcommand{\mP}{\mathcal{P}}

\newcommand{\bx}{{x}}

\newcommand{\by}{{y}}

\newcommand{\bw}{{w}}
\newcommand{\bt}{{\theta}}
\newcommand{\bxi}{{\xi}}

\newcommand{\bu}{{u}}

\newcommand{\RR}{\mathbb{R}}

\usepackage{calc}
\usepackage{accents}

\usepackage{booktabs,bm,float}
\usepackage{tcolorbox}

\usepackage{hyperref}
\usepackage[nameinlink]{cleveref}

\newtheorem{remark}{Remark}

\crefname{equation}{eq.}{eqns.}
\Crefname{equation}{Eq.}{Eqns.}

\newcommand{\revise}[1]{\textcolor{blue}{#1}}
\usepackage{subcaption}

\usepackage{makecell}
\begin{document}
\newtheorem{theorem}{Theorem}
\newtheorem{lemma}{Lemma}
\begin{frontmatter}

\title{Adversarial Numerical Analysis for Inverse Problems\tnoteref{t1}}

\author[label1]{Kailai Xu}
\ead{kailaix@stanford.edu}
\address[label1]{Institute for Computational and Mathematical Engineering, Stanford University, Stanford, CA, 94305}
 
\author[label1,label2]{Eric Darve}
\ead{darve@stanford.edu}
\address[label2]{Mechanical Engineering, Stanford University, Stanford, CA, 94305}

%
%
 
\begin{abstract}
Many scientific and engineering applications are formulated as inverse problems associated with stochastic models. In such cases the unknown quantities are distributions. The applicability of traditional methods is limited because of their demanding assumptions or prohibitive computational consumptions; for example, maximum likelihood methods require closed-form density functions, and Markov Chain Monte Carlo needs a large number of simulations. We introduce adversarial numerical analysis, which estimates the unknown distributions by minimizing the discrepancy of statistical properties between observed random process and simulated random process. The discrepancy metric is computed with a discriminative neural network. We demonstrated numerically that the proposed methods can estimate the underlying parameters and learn complicated unknown distributions. 
\end{abstract}

\begin{keyword}
Adversarial Training \sep Neural Networks \sep Automatic Differentiation
\end{keyword}

\end{frontmatter}

\section{Introduction}\label{sect:intro}

Almost all model-based problems in data analytics and scientific computing can be classified into two categories: forward problems and inverse problems. The forward problems assume that models and parameters---such as boundary conditions in partial differential equations---are explicitly given and ask for model outputs. However in the inverse problem, part of the parameters in the models are unknown but we observe (partial) outputs. The unknown parameters are calibrated based on stochastic models and observations. Among inverse problems, those associated with stochastic models are particularly interesting since stochasticity is indispensable for modeling real-world phenomenon. In the following we discuss several examples for inverse modeling with stochasticity
\begin{itemize}
        \item In \textbf{uncertainty quantification}, we have a partial differential equation (PDE) model
    \begin{equation}\label{model:1}
  \begin{aligned}
        \mathcal{L}_{\theta} u(x) &= 0 & x \in \Omega\\
        \mathcal{B} u(x) &= 0 & x \in \partial \Omega
    \end{aligned}
\end{equation}
    where $\mathcal{L}_{\theta}$ is a differential operator parametrized by an unknown $\theta$, and $\mathcal{B}$ is the operator associated with the boundary condition. Due to the stochastic nature of the physical process, $\theta$ is not deterministic and subject to an unknown distribution $\mathcal{P}$. Therefore, the resultant solution $u$ is a random variable. Given multiple observations $u_i$ for $u$, how can we reconstruct the unknown distribution $\mathcal{P}$? 
    \item In the \textbf{stochastic heat equation}
    \begin{equation}\label{model:2}
        \partial_t u(t, x) =  \nabla \cdot (\theta\nabla u(t, x)) + \xi, (t, x) \in \RR_+ \times \RR^d
    \end{equation}
    where $\xi$ is a centered Gaussian process such that $\mathbb{E}\xi(s,x)\xi(t,y) = \delta(t-s)\delta(x-y)$ ($\delta$ is the Dirac delta function). We have observed $u(t_i, x_i)$ at multiple times and locations, $i=1$, $2$, $\ldots$, $n$ and want to estimate the unknown conductivity coefficient $\theta\in\RR^{d\times d}$ from those observed data.
    \item In \textbf{finance}, the stochastic differentiation equations can be used for modeling interest rates. For example, the square root process given by
    \begin{equation}\label{model:3}
        dr_t = \alpha(\mu-r_t)dt + \sqrt{r_t}\sigma dW_t
    \end{equation}
    where $r_t$ is the interest rate, $\{W_t,t\geq 0\}$ is a standard Brownian motion and $\theta = (\alpha, \mu,\sigma)$ are model parameters. To describe the dynamics of the interest rates, one can calibrate the parameters $\theta$ based on discrete observations $r_{t_1}$, $r_{t_2}$, $\ldots$, $r_{t_n}$ at time $t=t_1$, $t_2$, $\ldots$, $t_n$. 
\end{itemize}

However, despite wide applications of inverse problems associated with stochastic models, they are challenging to solve and remain largely unexplored in the literature (compared to their \textit{deterministic} model counterparts). When the distributions in the models are non-Gaussian, there are few analytical tools available and therefore we must resort to numerical analysis. Traditionally, numerical methods have been mostly developed for inverse problems associated with deterministic models or Gaussian-based stochastic models. There are only a few approaches for non-Gaussian cases and their applicability is limited (see the following text).

Before we review traditional methods and describe our approach, we first formulate \Cref{model:1,model:2,model:3} by a unified mathematical model. Assume that $w$ is sampled from a known stochastic process and $\theta$ is a random variable with an unknown distribution. The governing equation for the system is 
\begin{equation}\label{equ:xfw}
    x = F(w, \theta)
\end{equation}
The output, $x$, will also be a random variable or process. Our task is to estimate the distribution of $\theta$, $\mathcal{P}$, from observations $\{x_i\}_{i=1}^n$. Even though the random process from which $w$ is sampled is known, $w$ itself is not observable; this is common in many applications. \Cref{tab:wtx} discusses the corresponding $w$, $\theta$ and $x$ in previous examples.

\begin{table}[hbt]
  \begin{tabular}{l|cccc}
  \toprule
  Example  & $w$& $\theta$ & $\mathcal{P}$& $x$ \\\midrule
   Uncertainty Quantification & -- & $\theta$ & Unspecified &$u(x)$ \\
   Stochastic Heat Equation & $\xi$ & conductivity $\theta$ & $\delta_{\theta}$ & $u(t,x)$ \\
   Stochastic Differential Equation & $W_t$ & $(\alpha, \mu, \sigma)$ & $\delta_{\alpha}\delta_{\mu}\delta_{\sigma}$ & $r_t$\\\bottomrule
  \end{tabular}
  \caption{The known random process $w$, unknown random variable $\theta$ and the corresponding distribution $\mathcal{P}$, and the output $x$. In the probability column, ``Unspecified'' means the form of the probability distribution is not specified, and $\delta_{\theta}$ means the Dirac delta distribution.}
  \label{tab:wtx}
\end{table}

We now have a brief review of traditional methods for inverse modeling of \Cref{equ:xfw} and below is their challenges and limitations. We will discuss two cases separately
\begin{enumerate}
    \item $\mathcal{P}=\delta_{\theta^*}(\theta)$ is a Dirac Delta distribution for an unknown $\theta^*$.
    \item The form of the distribution $\mathcal{P}$ of $\bt$ is not given; in this case, $\mathcal{P}$ can be approximated by a parametrized functional form and hence the problem is reduced to the last case.
\end{enumerate}
When $\mathcal{P}=\delta_{\theta^*}(\theta)$ is a Dirac Delta distribution for an unknown $\theta^*$, the inverse problem for \cref{equ:xfw} is reduced to the parameter estimation problem. The traditional starting point of estimating $\bt$ is the maximum likelihood estimator, which maximizes the log-likelihood function 
\begin{equation}
    \max_{\bt} \; l_{\bt}(x_1, x_2, \ldots, x_N) := \sum_{i=1}^N \log p(x_i | \bt)
\end{equation}
where $\{x_i\}_{i=1}^N$ are observed samples from the random process $x$. The maximum likelihood method~(MLE) is provably asymptotically efficient, that is, consistent and asymptotically normal with variance equal to the Cram\'er-Rao lower bound under certain conditions; typical conditions are given in \cite{daniels1961asymptotic,cramer1999mathematical,le1956asymptotic}. However, MLE requires computing $p(x_i|\bt)$. When $F$ is a very complicated model, it is difficult to obtain the analytical form of $p(x_i|\bt)$ or compute it numerically. 

In the general case where the form of the distribution $\mathcal{P}$ of $\bt$ is not given, we usually approximate $\mathcal{P}$ using various functional forms. Three approaches are popular in the literature
\begin{enumerate}
    \item Approximating $\mathcal{P}$ by a known distribution with tunable parameters. For example, Gaussian-mixture distributions are used in variational inference for approximating the unknown distribution $\mathcal{P}$~\cite{blei2017variational}. Consequently, the problem is reduced to estimating the covariances $\bm{\Sigma}$ and the means $\bm{\mu}$ of the mixed Gaussian distributions. The surrogate mathematical model can be formulated as 
\begin{equation}
    x = F(w, \tilde{\theta}(\bm{\Sigma}, \bm{\mu}))
\end{equation}
where $\tilde{\theta}(\bm{\Sigma}, \bm{\mu})$ is the Gaussian-mixture distribution and is used to approximate $\bt$. Even though $w$ is known, the log probability $p(x|\bm{\Sigma}, \bm{\mu})$ is usually analytically intractable due to the complicated nature of $F$, and therefore presents a major challenge for maximum likelihood methods for estimating $\bm{\Sigma}$, $\bm{\mu}$.
\item Bayes inference method~\cite{dashti2016bayesian}. In this approach, we update our estimation of $\bt$ according given each new observation $x_i$ to Bayes rule
\begin{equation}
    p(\theta|x_i) = \frac{p(x_i|\theta) p(\theta)}{p(x_i)}
\end{equation}
The Bayes inference method requires us to specify a prior $p(\bt)$, which can be difficult in the absence of a physical basis or a plausible scientific model; besides, it also requires calculating $p(x_i|\theta)$, which can be difficult to compute numerically. 
\item Markov Chain Monte Carlo~(MCMC). MCMC can create samples from $\mathcal{P}$ by constructing a Markov chain. It is a sample-based method and notoriously slow to converge, and sometimes suffers from non-convergence~\cite{van2018simple}.
\end{enumerate}


\revise{
\begin{table}[hbt]
  \begin{tabular}{l|cccc}
  \toprule
   Method &\makecell{Requires\\Sampling}  & \makecell{Requires\\Log Likelihood} &\makecell{Has Explicit\\Loss Functions} & \makecell{Parametrizes\\$\mathcal{P}$}  \\
    \midrule
  Gaussian-Mixture Approximation  & \xmark & \cmark & \cmark & \cmark  \\
  Bayes Inference  & \xmark  & \cmark & \cmark & \xmark \\
  MCMC  & \cmark & \xmark & \xmark & \xmark \\
  ANA &\cmark &\xmark & \cmark & \cmark \\ \bottomrule
  \end{tabular}
  \caption{Comparison of ANA and traditional methods.}
  \label{tab:compareANA}
\end{table}
}

Our approach, adversarial numerical analysis (ANA), is based on minimizing the discrepancy between the actual random process $x=F(w, \theta)$ and estimated random process $\tilde x=F(w, \tilde\theta)$. Here $\tilde \theta$ is our current estimation for the true parameter $\theta$ and is updated as the minimization proceeds. The key is that we can formulate the discrepancy between  statistical properties of $x$ and $\tilde x$  (called discriminator loss) with a neural network. The idea is borrowed from generative neural nets (GAN), where a neural network $D_{\xi}$, parametrized by $\xi$, tries to discriminate actual random variable $x$ and the estimated random variable $\tilde x$. Meanwhile, we update $\tilde \theta$ to generate samples $\tilde x$ that are indistinguishable from $x$. By updating the neural network and $\tilde\theta$ simultaneously, at equilibrium, we obtain a $\tilde \theta$ such that the neural network cannot discriminate $x$ and $\tilde x$, even though we keep updating $\xi$ to reinforce the neural network's discrimination ability. This indicates that the distribution $\tilde x$ is indistinguishable from $x$, under the probability metric determined by the neural network. 

How can the discriminator neural network measure the discrepancy between two probability distributions? The answer is that the choice of the loss functions $L^D(x, \tilde x; \xi)$ for updating $\xi$ (assuming $\tilde\theta$ is fixed), and the loss functions $L^F(\tilde x)=L^F(F(w, \tilde \theta))$ for updating $\tilde\theta$ (assuming $\xi$ is fixed) together determines the probability metric. For example, by properly choosing $L^D$ and $L^F$ (see \Cref{equ:klganloss} for such an example), we actually obtain the maximum likelihood estimator for $\theta$ (if $\theta$ is subject to a Dirac delta distribution) at equilibrium and the corresponding probability distance metric is Kullback-Leibler divergence. Even though we have not formulated it explicitly; more importantly, we do not have to compute log likelihood function analytically. It has been realized that some of the loss functions in GANs result in optimal estimators for specific probabilistic metrics for $D$. For example,
\begin{itemize}
    \item Vanilla GANs correspond to the Jensen-Shannon divergence. 
    \item Wasserstein GANs to the Wasserstein distance.
    \item Kullback-Leibler GANs to the KL divergence.
\end{itemize}

To approximate the unknown distribution $\mathcal{P}$ of $\theta$, we can parametrize the probability distribution with certain function forms. In the case when $\mP$ is \textbf{low-dimensional}, conventional functional forms such as linear basis functions, radial basis functions, and polynomial basis functions can be used to approximate the density function $\mP$. In this case, the problem is reduced to estimating some coefficients. 

\begin{figure}[hbtp]
  \includegraphics[width=0.8\textwidth]{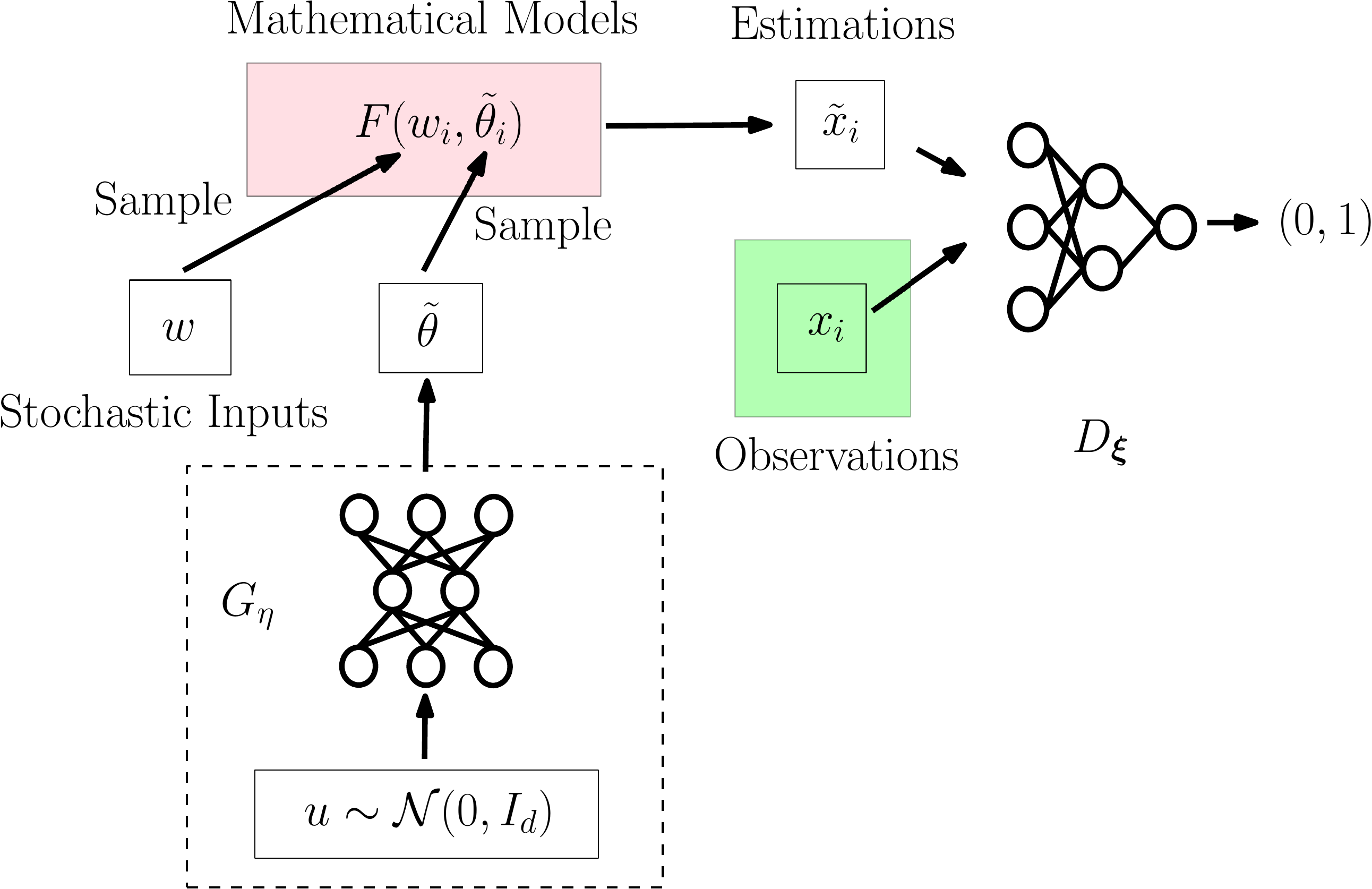}
  \caption{Graphical illustration of adversarial numerical analysis. $\{\bx_i\}$ are observations, $w_i$, $\tilde\theta_i$ are realizations of $w$ and $\tilde\theta$ respectively ($\tilde\theta=G_{\eta}(u)$ is an approximation to the true random variable $\theta$). $D_{\bxi}$ is the discriminator. We parametrize $\mP$ with a neural network $G_{\eta}$ and the problem is reduced to estimating $\eta$ (the dashed block). With estimated $\eta$, $G_{\eta}(u)$ is our approximation to the unknown random variable $\theta$. When $\mathcal{P}$ is a Dirac Delta distribution, the dashed block is not needed.}
  \label{fig:ana}
\end{figure}

However, when $\mP$ is a \textbf{high-dimensional} multi-variate distribution, those methods can be too computationally expensive. It has been found that neural networks are very good candidates to transform simple distributions to very complex ones. In GAN applications, such neural networks are called \textit{generative neural networks}. We adopt this approach in ANA and use a neural network $D_{\xi}:\RR^{d'}\rightarrow \RR^{d}$ to represent $d$-dimensional unknown distributions $\mP$ ($d'$ and $d$ are not necessarily equal). 
The generative neural net $G_{\eta}:\RR^{d'}\rightarrow \RR^d$ transforms a random variable $u$ whose distribution is easy to sample from, e.g., a multivariate Gaussian distribution $\mathcal{N}(0, I_{d'})$, to another random variable $G_{\eta}(u)$. We approximate $\theta$ with  $G_{\eta}(u)$ by tuning the parameters $\eta$. \Cref{fig:ana} illustrates the model formulation of ANA and \Cref{tab:compareANA} compares the difference between ANA and traditional approaches.

Now we discuss how ANA can be implemented. The idea is adversarial training, where we simultaneously (1) update $\tilde \theta$ such that the discriminator $D_{\xi}$ cannot distinguish $x$ and $\tilde x$; (2) update $\xi$ so such the discriminator is better at distinguishing $x$ and $\tilde x$. At equilibrium, the generated outputs $\tilde x$ will be indistinguishable from real random processes $x$ in the sense of small probability discrepancy, determined by $L^D$ and $L^F$. To be more concrete, assume that we have multiple observations $\{x_i\}$, the algorithm works as follows
\begin{enumerate}
    \item Generate realizations of $w$: $\bw_i$, $i=1$, $2$, $\ldots$, $N$.
    \item Generate realizations of $u\sim \mathcal{N}(0, I_{d'})$: $u_i$, $i=1$, $2$, $\ldots$, $N$.
    \item Compute $\tilde x_i = F(w_i, G_{\eta}(u_i))$.
    \item Compute the discrepancy between predictions $\{\tilde x_i\}$ and observations $\{ x_i \}$ in terms of probability distribution metrics: $d=L^D(\{ x_i \}, \{\tilde x_i\}; \xi)$, $i=1$, $2$, $\ldots$, $N$. 
    \item If $d$ is smaller than a predetermined threshold, stop and $G_{\eta}(u)$ is the approximation to $\theta$; otherwise, update $\eta$ according to the gradient
    \begin{equation}\label{equ:dtheta}
        \nabla_{\eta} L^F(\{F(w_i, G_{\eta}(u_i))\}_i)
    \end{equation}
    and repeat Step 1--4.
\end{enumerate}
Note we have abused the notation $L^D$, $L^F$ for both random variables/processes and discrete samples. The discrete version can be viewed as the numerical approximation to the random variables/processes counterparts. The explicit expression is given in \Cref{sect:choice}.

For our code implementation, we use automatic differentiation~\cite{baydin2018automatic} for computing \Cref{equ:dtheta}. This is key for implementing these methods quickly even for complex models. The model $F$ can be rather complicated: it can be a one-step simulation in the stochastic process, or a partial differential equation solver, or a sequence of coupled numerical procedures. For this purpose, we developed \texttt{ADCME}\footnote{Available at \url{https://github.com/kailaix/ADCME.jl}}, a \texttt{Julia} package, with \texttt{TensorFlow} and \texttt{PyTorch} backends, specially designed for engineering applications with automatic differentiation and parallel computing capabilities. It provides mathematically friendly syntax and supports hybrid programming with \texttt{C/C++} and \texttt{Julia}. It has the built-in optimizers and loss functions mentioned above. Those favorable features make it suitable for conducting adversarial numerical analysis for inverse modeling.

In the numerical examples, we first show that the method can calibrate an unknown parameter (a Dirac Delta distribution), and at the same time learn an unknown distribution. Then we demonstrate that the proposed algorithm is very effective in learning different complicated distributions, with the same neural network architectures for $D_{\xi}$ and $G_{\eta}$ respectively, for hidden physical parameters within a partial differential equation system.
Next, we benchmark different optimizers for ANA. We demonstrate numerically that the popular optimizer in engineering---\texttt{LBFGS}---perform equally well or even better in some cases compared to stochastic gradient descent optimizers. It converges faster and more smoothly than \texttt{ADAM} and \texttt{RMSProp}. We also discuss present numerical evidence for the convergence of the algorithm, which is consistent with our theoretical analysis. 
Throughout the numerical examples, we do not find significant differences in performance between the different loss functions.


 This paper is organized as follows: in \Cref{sect:ana} we introduce adversarial numerical analysis, and propose an optimization algorithm for the framework. In \Cref{sect:ca} we analyze the convergence of ANA based on Kullback-Leibler divergence, a special case of ANA for application in parameter calibration of CIR processes. In \Cref{sect:ne}, extensive numerical experiments are carried out, which include applications in uncertainty quantification, parameter inference, and option pricing. The numerical experiments also demonstrate our analysis in \Cref{sect:ca}. Finally, we conclude with a general discussion of the method in \Cref{sect:conc}.

\section{Adversarial Numerical Analysis}\label{sect:ana}

\subsection{Adversarial Numerical Analysis Description}\label{sect:def}

To illustrate the concept of adversarial numerical analysis, we consider a parameter inference problem in stochastic processes. The example is intended to be explanatory and ANA is applicable even when the closed-form density function is not available or the unknowns are probability distributions. 

Consider the parameter calibration of the mean in the Cox, Ingersoll and Ross process~(CIR process)~\cite{kladivko2007maximum}. The CIR process
\begin{equation}\label{equ:sim0}
    dr_t = \kappa(\theta-r_t)dt + \sqrt{r_t}\sigma dW_t
\end{equation}
is used for interest rate modeling; $(\kappa$, $\theta$, $\sigma)$ are model parameters. $\{W_t, t\geq 0\}$ is the standard Brownian motion; the interest rate $r_t$  moves in the direction of its mean $\theta$ at speed $\kappa$, $r_t\sigma^2$ is the diffusion function. Assume we are given a sample path $\mathbf{R} = \{R_1, R_2, \ldots, R_n\}$ with time interval $\Delta t$, and $\sigma$, $\kappa$ are known. The task is to estimate the mean $\theta$. Many methods exist for such a problem, such as maximum likelihood methods~\cite{kladivko2007maximum} and ordinary least square methods~\cite{anderson2003comparison}.

\begin{figure}[hbt]
  \includegraphics[width=0.8\textwidth]{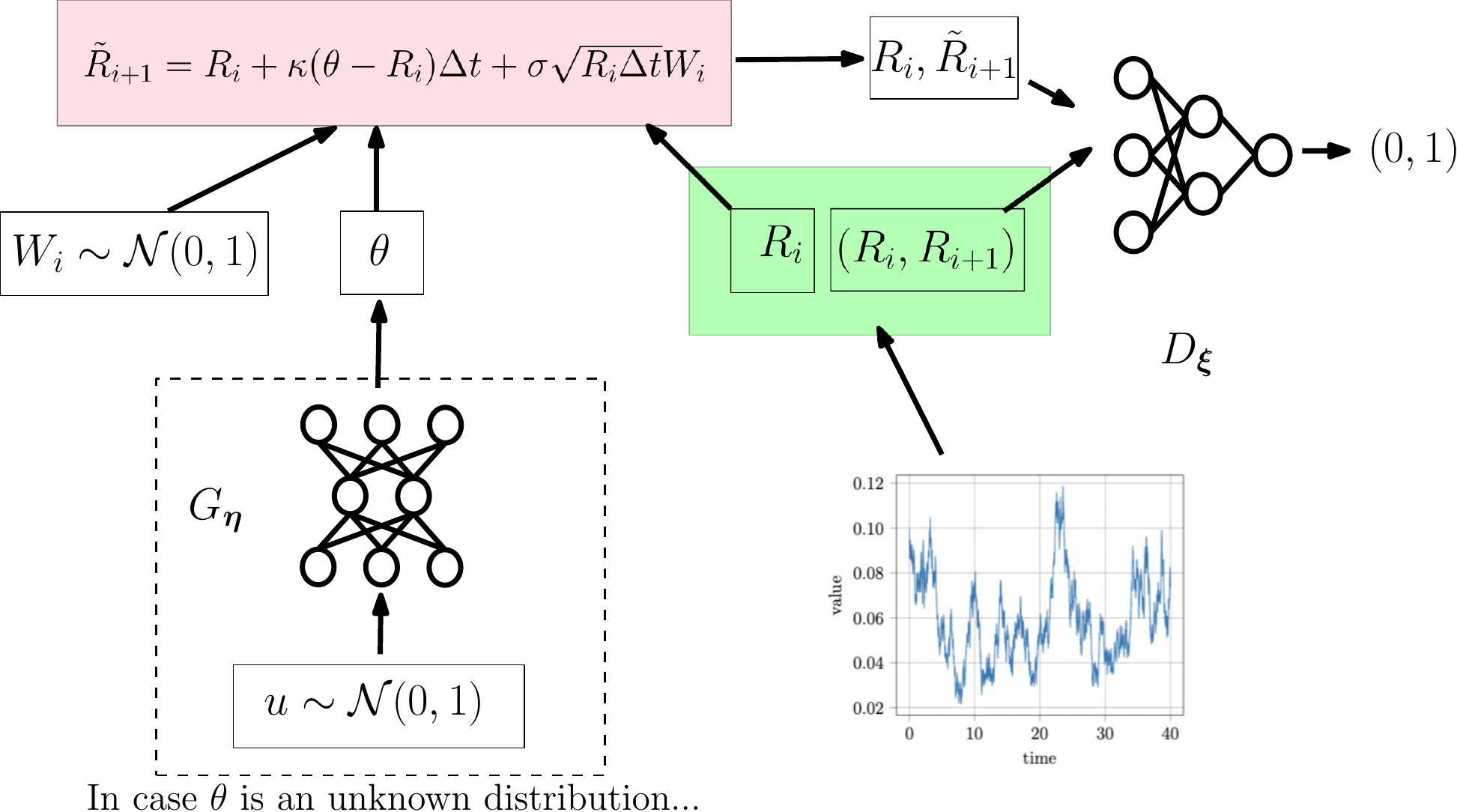}
  \caption{ANA pipeline for \Cref{equ:sim0}. The outputs are obtained from simulated samples $\tilde R_{i+1}$ from inputs $R_i$. $D_{\bxi}$ discriminates $\{(R_i, \tilde R_{i+1})\}$ and $\{(R_i, R_{i+1})\}$.}
\end{figure}

\ed{not clear why you have $R_i$ in the green box; you should separate the ICs for your ODE and the input to $D_\xi$.}

The idea of ANA is to match the distribution of simulated data and that of observed data. For this example, we define the observed data as $\{(R_i,R_{i+1})\}_{i=1}^n$. The numerical model corresponding to the continuous one \Cref{equ:sim0} is the Euler Maruyama scheme
\begin{equation}
    \tilde\by_i = F((x_i, W_i), \theta) = \left( x_i,\ x_i + \kappa (\theta-x_i) \Delta t + \sigma \sqrt{x_i \Delta t}W_i \right)
\end{equation}
with stochastic input $W_i \sim \mathcal{N}(0,1)$ and $x_i$ is randomly drawn from the sample path. A discriminative neural network $D_{\bxi}: \RR^2\rightarrow (0,1)$ is used to measure the discrepancy between $\{\by_i\}$ and $\{\tilde\by_i\}$. The goal of the neural network is to output 1 for $\by_i$ and $0$ for $\tilde \by_i$
\begin{equation}
    D_{\bxi}(\by_i) \approx 1\qquad D_{\bxi}(\tilde\by_i)\approx 0
\end{equation}
under this assumption; one choice of the model loss function is~(recall that $\tilde\by_i$ depends on $\theta$)
\begin{equation}\label{equ:modelloss}
    L^F(\{ \tilde\by_i \}) = -\frac{1}{n}\sum_{i=1}^n \log D_{\bxi}(\tilde\by_i)
\end{equation}

when we drive $D_{\bxi}(\tilde\by_i)$ to 1 by minimizing $L^F$, and therefore making the simulated outputs $\tilde\by_i$ seem more similar to the observed ones. 

By parametrizing $\mP$ with the neural network $D_{\bxi}$, the mathematical model for the stochastic process is converted to 
\begin{equation}\label{equ:reduced}
    \tilde\by_i = F((x_i, W_i), D_{\bxi}(u)) 
\end{equation}
here $u \sim \mathcal{N}(0,1)$, but it can also be other probability distributions such as uniform distributions, as long as we are able to draw samples. The optimization algorithm follows \Cref{algo:ana}. We can apply \Cref{algo:ana} to the new model \Cref{equ:reduced}. 

\begin{algorithm}[htpb]
\caption{{Optimization algorithm for ANA. The outer loop of the algorithm involves updating $\bt$ and an inner loop, where the discrimination neural network is updated for $k$ times. For \texttt{LBFGS}, $k=1$ is used.}}
\label{algo:ana}
\begin{algorithmic}[1]
\FOR{$i=1,2,3,\ldots$}
\FOR{$j=1,2,\ldots,k$}
\STATE Generate $n$ samples $\{\bw_1,\bw_2,\ldots,\bw_n\}$ from $\mathcal{Q}$
\STATE Compute $\tilde\by_i\gets F(\bx_i, \bw_i; \bt)$, $i=1$, $2$, $\ldots$, $n$
\STATE Update the discrimination neural network parameters $\bxi$ by the gradients
$$    \nabla_{\bxi} L^D(\{\by_i\},\{\tilde\by_i\}) $$
\ENDFOR
\STATE Generate $n$ samples $\{\bw_1,\bw_2,\ldots,\bw_n\}$ from $\mathcal{Q}$
\STATE Compute $\tilde\by_i\gets F(\bx_i, bw_i; \bt)$, $i=1$, $2$, $\ldots$, $n$
\STATE Compute the gradients of $L^F$ using adjoint state methods~(or automatic differentiation)
$$    \mathbf{g} = \nabla_{\bt} L^F(\{\tilde\by_i\}) $$
\STATE Update $\bt$ with $\mathbf{g}$
\ENDFOR 
\end{algorithmic}    
\end{algorithm}
In our algorithm, we have used multiple optimizers for updating $\bxi$ and $\bt$, such as gradient descent and \texttt{LBFGS}.

We remark that technically the method can also be used to estimate unknown parameters and distributions at the same time. Besides, we have not discussed the well-posedness and the conditioning of the problems yet but we will acknowledge its importance through the analysis of a model problem in \Cref{sect:ca}.

\subsection{Choice of Loss Functions for ANA}\label{sect:choice}

The discriminator neural network $D_{\bxi}$ shares similar features as that in GANs. In this subsection, we propose three sets of loss functions, which all come from different GAN algorithms.

Recall that $D_{\bxi}$ is the discrimination neural network, $\by_i$ and $\tilde\by_i$ are observed and simulated data respectively. The major difference between the following GANs is their loss functions. 
\begin{itemize}
    \item \textit{Vanilla GAN}~\cite{goodfellow2014generative}. The output of $D_{\bxi}$ is a number in the range $(0,1)$, and the loss functions are 
    \begin{align}
        L^F(\{\tilde\by_i\}) =& -\frac{1}{n}\sum_{i=1}^n \log D_{\bxi}(\tilde\by_i)\\
        L^D(\{\by_i\}, \{\tilde\by_i\}; \xi) =& -\frac{1}{n}\sum_{i=1}^n\left(\log(D_{\xi}(\by_i))  + \log (1-D_{\xi}(\tilde \by_i))\right)
    \end{align} 

    At equilibrium, $L^F = \log 2$, $L^D = \log 4$.
    \item \textit{Wasserstein GAN}~\cite{arjovsky2017wasserstein}. The last layer of the discrimination network is usually a linear layer, and the loss functions are 
    \begin{align}
        L^F(\{\tilde\by_i\}) =& \frac{1}{n}\sum_{i=1}^n D_{\xi}(\tilde \by_i) \\
        L^D(\{\by_i\}, \{\tilde\by_i\}; \xi) =& -\frac{1}{n}\sum_{i=1}^n D_{\xi}(\by_i)
    \end{align}
    The discrimination neural network $D_{\bxi}$ is usually restricted to those whose norm of weights are within $[-c,c]$, $c>0$.  The norm of the weights in $D_{\bxi}$ are usually restricted to $[-c,c]$, $c>0$, so that $D_{\bxi}$ belongs to the class of Lipschitz functions. In this case, the output of the neural network is not necessarily $(0,1)$. At equilibrium, $L^F = 0$, $L^D = 0$. 

    \item \textit{Kullback–Leibler~(KL) GAN}~\cite{nowozin2016f}. The output of $D_{\xi}$ is a number in the range $(0,1)$, and the loss functions are 
    \begin{equation}\label{equ:klganloss}
    \begin{aligned}
            L^F(\{\tilde\by_i\}) =& \frac{1}{n}\sum_{i=1}^n \log\frac{ (1-D_{\xi}(\tilde \by_i))}{D_{\xi}(\tilde \by_i)} \\
            L^D(\{\by_i\}, \{\tilde\by_i\}; \xi) =& -\frac{1}{n}\sum_{i=1}^n \left( \log D_{\xi}(\tilde \by_i) + \log(1-D_{\xi}(\by_i)) \right)
        \end{aligned} 
    \end{equation}
    At equilibrium, $L^F = 0$, $L^D = \log 4$. 
\end{itemize}

Other variations of loss functions exist. A comprehensive comparison of the loss functions is out of scope but we refer readers to \cite{nowozin2016f} for an excellent discussion. 

The choice of optimizers is another concern for training in ANA. We attempt to compare the \texttt{LBFGS} optimizer, which is extensively used in engineering, with other popular optimizers in machine learning. In our experiment, we mainly use three kinds of optimizers
\begin{itemize}
    \item \texttt{ADAM}~\cite{kingma2014adam}. Adam is an algorithm for first-order optimization of stochastic objective functions based on adaptive estimates of lower-order moments. We apply 5 updates to $\bt$ and 1 update to $D_{\bxi}$ in each iteration. 
    \item \texttt{RMSProp}~\cite{lectures70:online}. \texttt{RMSProp} divides the gradient by a running average of its recent magnitude. Usually it divides the learning rate by an exponentially decaying average of squared gradients.  We apply 5 updates to $\bt$ and 1 update to $D_{\bxi}$ in each iteration. 
    \item \texttt{LBFGS}~\cite{skajaa2010limited}. The limited Broyden–Fletcher–Goldfarb–Shanno algorithm is a powerful approach for finding a local minimum of a nonconvex objective function. The \texttt{LBFGS} approximates the Hessian based on the most recent gradients. In our examples, \texttt{LBFGS} is used for optimizing $L^F$ and one \texttt{ADAM} update is applied to $D_{\bxi}$ in each iteration. 
\end{itemize}

\section{Convergence Analysis}\label{sect:ca}
\subsection{KL GAN Produces Maximum Likelihood Estimators}\label{sect:mle}

{Thanks to the connection of ANA with GANs discussed in \Cref{sect:intro}, many convergence analysis results for GANs are also applicable to ANA.} There are many variants of GANs. From the perspective of divergence minimization~\cite{jolicoeur2018gans}, variants are proposed to minimize Kullback-Leibler distance~(KL), the Jensen-Shannon distance, or the Wasserstein distance. For example, the vanilla GAN proposed in \cite{goodfellow2014generative} is equivalent to minimizing the cross-entropy~\cite{goodfellow2016nips}.
One particularly interesting case is minimizing the KL-divergence.  It was shown \cite{goodfellow2016nips} that it is equivalent to maximizing the log-likelihood function
\begin{equation}\label{equ:kldiv}
    \arg\min_{ \eta} D_{KL}(\{x_i\}_{i=1}^n || F( w, G_{\eta}(u) )) = \arg\max_{\eta} \log \sum_{i=1}^n l_{\eta}(x_i)
\end{equation}
where $l_{\eta}(x)$ is the log likelihood function and $\{x_i\}_{i=1}^n$ is the discrete distribution of the observations $x$. The optimal $\eta$ is exactly the maximum likelihood estimator. \cite{goodfellow2014generative} showed that using \cref{equ:klganloss} is equivalent to minimizing \cref{equ:kldiv} in expectation, under the assumption that the discrimination neural network is optimal.

    It is known that under some regularity conditions on the family of distributions, the maximum likelihood estimator $\hat{\eta}_n$ is consistent, i.e., if ${\eta}^*$ is the true parameter, we have $\hat{\eta}_n\rightarrow {\eta}^*$, $n\rightarrow\infty$~\cite{myung2003tutorial}. In addition, we have the asymptotic normality for maximum likelihood estimator~\cite{ly2017tutorial}
\begin{equation}
    \sqrt{n}(\hat {\eta}_n - {\eta}^*) \rightarrow \mathcal{N}\left(0, \frac{1}{I({\eta}^*)} \right)
\end{equation}
where $I({\eta})$ is the Fisher information
\begin{equation}
    I({\eta}) = - \mathbb{E}_{\bx\sim \mathcal{U}_\mathcal{X}}\left(\frac{\partial^2}{\partial{\eta}^2} \log l_{\eta}(\bx)\right)
\end{equation}

\begin{remark}
    The discussion above holds true if $\theta$ is subject to a Dirac delta distribution. In this case, the minimization problem becomes 
    \begin{equation}\label{equ:kldiv}
    \arg\min_{ \theta} D_{KL}(\{x_i\}_{i=1}^n || F( w, \theta )) = \arg\max_{\theta} \log \sum_{i=1}^n l_{\theta}(x_i)
\end{equation}
and let the maximum likelihood estimator for $n$ samples be $\hat\theta_n$, and the exact parameter is $\theta^*$, then we have
\begin{equation}
    \sqrt{n}(\hat \theta_n - \theta^*) \rightarrow \mathcal{N}\left( 0, \frac{1}{I(\theta^*)} \right)
\end{equation}
\end{remark}

\subsection{Theoretical Analysis for the CIR Example with KL Divergence}\label{sect:cireg}

In this section, we analyze the convergence of ANA for estimating $\theta$ and $\kappa$ in \Cref{equ:sim0}. The analysis for other parameters is similar. The analysis also gives us some insights into when the mathematical problem is proposed and thus guides the design of algorithms.

In the last section \Cref{sect:mle} we showed that if we use KL GAN, under the assumption that the discrimination neural network is optimal, $\kappa$, $\theta$ will converge to the maximum likelihood estimator in expectation. We assume that all those assumptions are satisfied, and therefore we can focus on analyzing the convergence of the maximum likelihood estimators for $\theta$ and $\kappa$. Note the choice of KL GAN is merely for analysis purposes and reconstructing the maximum likelihood estimator from ANA is not the final goal. But since the maximum likelihood estimator is well studied the analysis based on this particular choice helps us understand the theoretical aspect of ANA better. 
\subsubsection{Convergence for the $\tau$ estimator}
We revisit the CIR example in \Cref{sect:def}. Recall in the CIR example, $y = F(w, \theta)$, samples of the known stochastic process $w$ have the form $W_i \sim \mathcal{N}(0,1)$ and the output $(x_i, y_i)$ are two consecutive sample with time interval $\Delta t$. $x_i$ and $y_i$ are related by
\begin{equation}
    y_i = x_i + \kappa(\tau-x_i)\Delta t + \sigma\sqrt{x_i\Delta t}W_i
\end{equation}

 The probability density function of $y_i$ given the observation $x_i$ is given by 
\begin{equation}
    p(y|x) = {1 \over {\sqrt {2\pi \sigma^2 x\Delta t} }}\exp \left( { - {1 \over 2}{{\left( {{{y - x - \kappa (\tau  - x)\Delta t} \over {\sigma \sqrt {x\Delta t} }}} \right)}^2}} \right)
\end{equation}

Assuming that $x$ has the probability density function $f(x)$, then the log likelihood function for the joint distribution $(x,y)$ is 
\begin{equation}
    l_\tau(x, y) = - \frac{1}{2}{\left( {\frac{{y - x - \kappa (\tau  - x)\Delta t}}{{\sigma \sqrt {x\Delta t} }}} \right)^2} + \log {f}(x) - \log \left( {\sigma \sqrt {2\pi x\Delta t} } \right)
\end{equation}
Therefore, assume we have observations $(x_i, y_i)$, $i=1$, $2$, $\ldots$, $n$, which are consecutive samples with time interval $\Delta t$, the empirical log likelihood function is
\begin{equation}
    \sum_{i=1}^nl_\tau(x_i,y_i) = - {1 \over 2}\sum\limits_{i = 1}^n {{{\left( {{{{y_i} - {x_i} - \kappa (\tau  - {x_i})\Delta t} \over {\sigma \sqrt {{x_i}\Delta t} }}} \right)}^2}}  + \sum\limits_{i = 1}^n \log {{f}({x_i})}  - \sum_{i=1}^n\log \left( {\sigma \sqrt {2\pi x_i\Delta t} } \right)
\end{equation}
the maximum likelihood estimator can be computed
\begin{equation}\label{equ:ml}
    \sum_{i=1}^nl'_\tau(x_i,y_i) = 0\Rightarrow \hat\tau_n = {{{1 \over n}\sum\limits_{i = 1}^n {{{{y_i}} \over {{x_i}}} + \kappa \Delta t - 1} } \over {\left( {{1 \over n}\sum\limits_{i = 1}^n {{1 \over {{x_i}}}} } \right)\kappa \Delta t}}
\end{equation}

According to the discussion, if we use KL divergence as the discrepancy measure of the real sample-path and generated ones, $\hat \tau$ at the Nash equilibrium will converge to $\hat \tau_n$ under mild assumptions. Let $\tau^*$ be the exact solution, we prove that as $\Delta t\rightarrow 0$, $n\rightarrow \infty$, $\hat\tau_n\rightarrow \tau^*$. We always assume convergence in the distribution for probability convergence.

\begin{theorem}
Assume that $\mathbb{E}\left(\frac{1}{x}\right) = X_{-1}\in (0, \infty)$, for sufficiently small $|\Delta t|$, we have
    \begin{equation}\label{equ:dt}
        \hat \tau_n \rightarrow \tau^*+\mathcal{O}(\Delta t) \qquad n \rightarrow \infty
    \end{equation} 
    in addition, the corresponding Fisher information  
    \begin{equation}
        I(\tau) = \frac{\kappa^2\Delta tX_{-1}}{\sigma^2    }
    \end{equation}
    and consequently 
    \begin{equation}
        \sqrt{n}(\hat \tau_n -\tau^*) \rightarrow \mathcal{N}\left( 0,\frac{\kappa^2\Delta t}{\sigma^2    X_{-1} } \right)
    \end{equation}
\end{theorem}

\begin{proof}
     
        Since $\mathbb{E}\left(\frac{1}{x}\right)<\infty$, by law of large numbers we have
    \begin{equation}\label{equ:1}
        \frac{1}{n}\sum_{i=1}^n\frac{1}{x_i}\rightarrow X_{-1}\qquad n \rightarrow \infty
    \end{equation}
    
    We first verify that  $\mathbb{E}\left(\frac{y}{x}\right)<\infty$ so that the law of large numbers hold. For CIR processes, we have~\cite{chou2006some}
    \begin{equation}
        {\Bbb E}\left[ {y|x} \right] = x{e^{ - \kappa \Delta t}} + \tau^* (1 - {e^{ - \kappa \Delta t}})
    \end{equation}
    thus 
    \begin{equation}
        {\Bbb E}\left[ \frac{y}{x} \right] = {\Bbb E}\left[{\Bbb E}\left[ \frac{y}{x}\Big|x \right]\right] = e^{-\kappa\Delta t} + \tau^* (1-e^{-\kappa \Delta t})X_{-1} <\infty
    \end{equation}
    therefore
    \begin{equation}\label{equ:2}
        \frac{1}{n}\sum_{i=1}^n\frac{y_i}{x_i}\rightarrow e^{-\kappa\Delta t} + \tau^* (1-e^{-\kappa \Delta t})X_{-1} \qquad n\rightarrow \infty
    \end{equation}
    Plug \cref{equ:1,equ:2} into \cref{equ:ml} we obtained using the continuous mapping theorem~\cite{Apr5pdf22:online}
    $$\hat \tau_n \rightarrow \frac{{{\tau ^*}(1 - {e^{ - \kappa \Delta t}})}}{{\kappa \Delta t}} + \frac{{{e^{ - \kappa \Delta t}} + \kappa \Delta t - 1}}{{{X_{ - 1}}\kappa \Delta t}} = {\tau ^*} + {\mathcal O}(\Delta t)$$
    \end{proof}

\begin{remark}
The constraint $\mathbb{E}\left(\frac{1}{x}\right) = X_{-1}<\infty$ is satisfied asymptotically.  Assume $2\kappa \tau > \sigma^2$, and $n$ is large enough so $x_i$ is subject to the asymptotic distribution of the CIR process, which is a gamma distribution  with density function 
    \begin{equation}\label{equ:hr}
        h(r) = \frac{w^\nu}{\Gamma(\nu)}r^{\nu-1} e^{-wr}, \quad 
        w = \frac{2\kappa}{\sigma^2}, \quad
        \nu = \frac{2\kappa \tau^*}{\sigma^2}
    \end{equation}
thus we have
    \begin{equation}\label{equ:hreq}
        \mathbb{E}\left(\frac{1}{x}\right) = \int_0^\infty \frac{1}{r}h(r)dr = \frac{1}{\tau^*}<\infty
    \end{equation}
\end{remark}

If $x$ is sampled from the asymptotic distribution (when $n$ is very large), the Fisher information is
\begin{equation}
    I(\tau) = \frac{\kappa^2\Delta t}{\sigma^2\tau^*}
\end{equation}
therefore we have
\begin{equation}\label{equ:fisher}
    \sqrt{n}(\hat \tau_n - \tau^*)\rightarrow \mathcal{N}\left(0, \frac{\sigma^2\tau^*}{\kappa^2\Delta t} \right)
\end{equation}

The condition \cref{equ:dt} and \cref{equ:fisher} indicates that for $\hat\tau_n\rightarrow \tau^*$ as $\Delta t\rightarrow 0$, $n\rightarrow \infty$, the following condition is sufficient
\begin{equation}\label{equ:condd}
    n\Delta t \rightarrow \infty,\ \Delta t \rightarrow 0
\end{equation}
for example, we let $\Delta t= \frac{1}{\sqrt{n}}$, i.e., $n = \frac{1}{\Delta t^2}$, then the condition \cref{equ:condd} is satisfied.

\subsubsection{Convergence for the $\kappa$ estimator}

We assume that $\tau$ and $\sigma$ are known and we want to estimate $\kappa$, whose true value is $\kappa^*$. 
\begin{theorem}\label{thm:div}
    Assume that $\mathbb{E}(x)=X_0\in(0,\infty)$, $\mathbb{E}\left(\frac{1}{x} \right)=X_{-1}<\infty$, for sufficiently small $|\Delta t|$, we have
    \begin{equation}
        \hat\kappa_n \rightarrow \kappa^* + \mathcal{O}(\Delta t) \qquad n\rightarrow \infty
    \end{equation}
    in addition, the corresponding Fisher information  
    \begin{equation}
        I(\kappa) = \frac{\Delta t}{\sigma^2    }(\tau^2 X_{-1} -2\tau + X_0)
    \end{equation}
    and consequently 
    \begin{equation}\label{equ:divkappa}
        \sqrt{n}(\hat \kappa_n -\kappa^*) \rightarrow \mathcal{N}\left( 0, \frac{\sigma^2}{\Delta t(\tau^2 X_{-1} -2\tau + X_0 )} \right)\qquad n\rightarrow \infty
    \end{equation}
\end{theorem}
\begin{proof}
The log-likelihood function for estimating $\kappa$ is 
\begin{equation}
    l_\kappa(x,y) = - \frac{1}{2}{\left( {\frac{{y - x - \kappa (\tau  - x)\Delta t}}{{\sigma \sqrt {x\Delta t} }}} \right)^2} + \log {f_1}(x) - \log \left( {\sigma \sqrt {2\pi x\Delta t} } \right)
\end{equation}
and therefore the empirical log likelihood function is
\begin{equation}
    \sum_{i=1}^nl_\kappa(x_i,y_i) = - {1 \over 2}\sum\limits_{i = 1}^n {{{\left( {{{{y_i} - {x_i} - \kappa (\tau  - {x_i})\Delta t} \over {\sigma \sqrt {{x_i}\Delta t} }}} \right)}^2}}  + \sum\limits_{i = 1}^n \log {{f_1}({x_i})}  - \sum_{i=1}^n\log \left( {\sigma \sqrt {2\pi x_i\Delta t} } \right)
\end{equation}
the maximum likelihood estimator can be computed \begin{equation}\label{equ:ml2}
    \sum\limits_{i = 1}^n {{{l'}_\kappa }} ({x_i},{y_i}) = 0 \Rightarrow {\hat \kappa _n} = \frac{1}{{\Delta t}}\frac{{\tau \frac{1}{n}\sum\limits_{i = 1}^n {\frac{{{y_i}}}{{{x_i}}} - \frac{1}{n}\sum\limits_{i = 1}^n {{y_i}}  - \tau  + \frac{1}{n}\sum\limits_{i = 1}^n {{x_i}} } }}{{{\tau ^2}\frac{1}{n}\sum\limits_{i = 1}^n {\frac{1}{{{x_i}}} - 2\tau  + \frac{1}{n}\sum\limits_{i = 1}^n {{x_i}} } }}
\end{equation}

    By the law of large numbers we have 
    \begin{equation}\label{equ:proof1}
        \frac{1}{n}\sum_{i=1}^n\frac{1}{x_i}\rightarrow X_{-1} \qquad \frac{1}{n}\sum_{i=1}^n{x_i}\rightarrow X_{0}
    \end{equation}
    in addition, we have
    \begin{equation}\label{equ:proof2}
        \frac{1}{n}\sum_{i=1}^n{y_i}\rightarrow \mathbb{E}[{\mathbb{E}[y|x]}] = X_0 e^{-\kappa^* \Delta t} + \tau (1-e^{-\kappa^* \Delta t})
    \end{equation}
    Plug \cref{equ:proof1,equ:proof2} into \cref{equ:ml2} and note $e^{-\kappa^*\Delta t} \rightarrow 1-\kappa^* \Delta t$, $\Delta t\rightarrow 0$, we have
    \begin{align}
        \hat \kappa_n  \rightarrow & \frac{1}{{\Delta t}}\frac{{\tau \left( {{e^{ - {\kappa ^*}\Delta t}} + \tau (1 - {e^{ - {\kappa ^*}\Delta t}}){X_{ - 1}}} \right) - \left( {{X_0}{e^{ - {\kappa ^*}\Delta t}} + \tau (1 - {e^{ - {\kappa ^*}\Delta t}})} \right) - \tau  + {X_0}}}{{{\tau ^2}{X_{ - 1}} - 2\tau  + {X_0}}}\\
        \rightarrow & \kappa^* + \mathcal{O}(\Delta t) \qquad \Delta t \rightarrow 0
    \end{align}
    
    Since
    \begin{equation}
        l''_\kappa(x) = -\frac{(\tau-x)^2\Delta t}{\sigma^2 x}
    \end{equation}
    we have
    \begin{equation}
        I(\kappa) = -\mathbb{E}(l''_\kappa(x)) = -\frac{\Delta t}{\sigma^2}\left( \tau^2 \mathbb{E}\left( \frac{1}{x} \right) - 2\tau + \mathbb{E}(x) \right) = \frac{\Delta t}{\sigma^2    }(\tau^2 X_{-1} -2\tau + X_0)
    \end{equation}
\end{proof}

In the case where $n$ is very large, $\frac{1}{n}\sum_{i=1}^n\frac{1}{x_i}$ and $\frac{1}{n}\sum_{i=1}^n{x_i}$ are approximately subject to the asymptotic distribution. Therefore, we have from \cref{equ:hreq}
\begin{equation}
    X_{-1} \approx \frac{1}{\tau}\qquad X_0 \approx \tau
\end{equation}

Consequently
\begin{equation}
    I(\kappa) \approx 0
\end{equation}
which indicates that the estimator has infinite variance. However, as long as we resample from the observations such that the new training data satisfies $\tau^2 X_{-1} -2\tau + X_0\neq 0$, we will obtain convergence result for $\kappa$ estimator. One such technique is given in \Cref{sect:num_cir}.
 
\begin{remark}
    Given a parameter $\tau$, let the family of distribution $P_\tau$ and $P_{\tau^*}$ be the true data distribution. The relationship between KL divergence and Fisher information can be described as follows for sufficiently small $|\tau-\tau^*|$ under mild assumptions
    \begin{equation}
        D_{KL}(P_{\tau^*} || P_\tau) \approx \frac{1}{2} (\tau-\tau^*)^T I(\tau^*) (\tau-\tau^*)
    \end{equation}
    thus if $I(\tau^*)\approx 0$, we can see that $D_{KL}(P_{\tau^*} || P_\tau)$ will have vanishing gradients and Hessians at $\tau=\tau^*$. This may lead to difficulty during optimization. 
\end{remark}

\section{Numerical Examples}\label{sect:ne}

In this section, we present four numerical examples. The mathematical models for those problems can all be formulated in the form of $x = F(\bw, \bt)$ and therefore ANA is applicable for those cases. The extensive numerical results also demonstrate the generality of the method. 

\subsection{Estimation and Uncertainty Quantification of Hidden Parameters in PDEs}\label{sect:n1}

In this example, the unknown parameter $\bt$ consists of an unknown distribution and an unknown parameter. The unknown distribution is parametrized by a neural network. We simultaneously estimate both the distribution and the parameter. 

Consider a Poisson equation with Dirichlet boundary condition
\begin{equation}\label{equ:poisson}
\begin{cases}
    -\nabla \cdot (a(x)\nabla u(x)) = 1 & x\in(0,1)\\
    u(0) = u(1) = 0 & \mbox{otherwise}
\end{cases}
\end{equation}
where
\begin{equation}
    a(x) = 1-0.9\exp\left( -\frac{(x-\mu)^2}{2\sigma^2} \right)
\end{equation}
Assume that $\mu$ and $\sigma$ are both unknown, and there is model form uncertainty in $\mu$, i.e., $\mu\sim \mP$ for some unknown probability distribution $\mP$. Our goal is that given many observations $\mathbf{u}_i$, $i=1$, $2$, $\ldots$, $N$, where $\mathbf{u}_i\in \mathbb{R}^n$ are observations $u(h)$, and each $\bu_i$ consists of $n=100$ nodal values of $u(x)$, i.e.,\ $u(2h)$, $\ldots$, $u(nh)$, $h=\frac{1}{n+1}$, we want to determine $\sigma$ and $\mP$.

Traditionally, this is usually done under the Bayesian framework where a priori information is imposed on $\mP$ and the transformation of probability in the forward model is calculated or approximated. However, this approach may be infeasible for complicated problems where the probabilistic relationship between data and parameter is analytically intractable, i.e., expressing $p(u(h), u(2h), \ldots, u(nh)| \ \mu, \sigma)$ in a closed-form. We present an example where ANA can be used for solving such problems. 

\begin{figure}[hbtp]
  \includegraphics[width=0.8\textwidth]{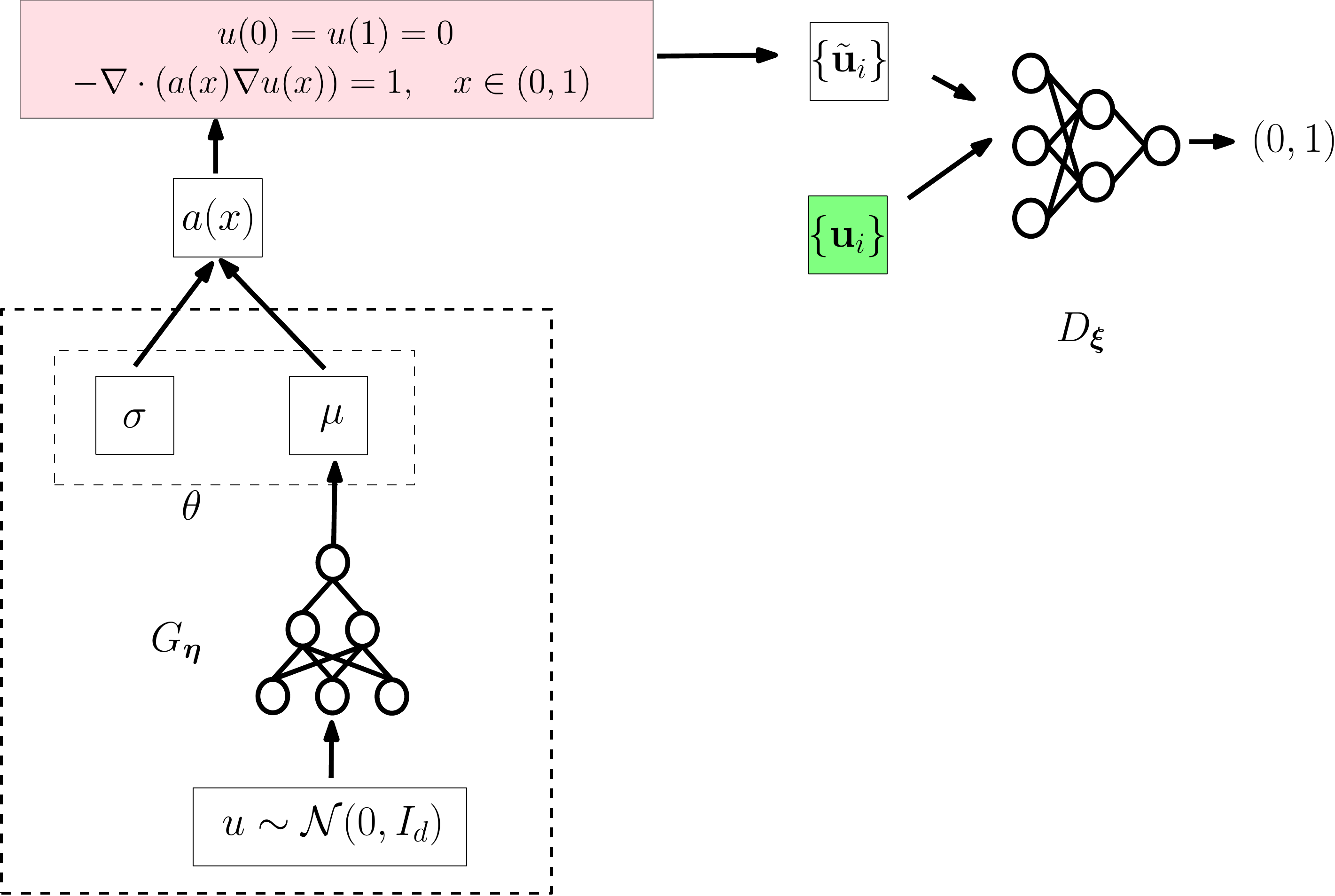}
  \caption{The ANA computation pipeline for \cref{equ:aaa}.}
  \label{fig:poissonpipe}
\end{figure}

To parametrize $\mP$, we consider a neural network $N:\mathbb{R}^d\rightarrow \mathbb{R}$, $d=10$; the inputs of neural network are samples from the uniform distribution $\mathcal{U}([-1,1]^d)$ while the output is a generated sample of $\mu$. To discretize  \cref{equ:poisson}, we consider the central difference scheme
\begin{equation}\label{equ:aaa}
  \begin{aligned}
    -a_{i-\frac{1}{2}} u_{i-1} + (a_{i-\frac{1}{2}}+a_{i+\frac{1}{2}})u_i - a_{i+\frac{1}{2}} u_{i+1} & = h^2,  & \quad i = 1, 2, \ldots, n\\
    u_0 = u_{n+1} & = 0
\end{aligned}
\end{equation}
where $a_{j\pm \frac{1}{2}} = a( h(j\pm \frac{1}{2}))$ which depends on $\mu$ and $\sigma$. The equation \cref{equ:aaa} leads to a tridiagonal system and can be solved using the Thomas algorithm~\cite{golub2012matrix}. The computation pipeline is shown in \cref{fig:poissonpipe}. The loss functions for $F_\theta$ and the discrimination neural network are chosen to be Wasserstein GAN losses. We use \texttt{RMSProp} with learning rate $10^{-4}$ and batch size 32. The true $\sigma^*=0.1$ and the true distribution for $\mu$ is 
\begin{equation}
    \mu^* \sim \mathcal{N}(0.3, 0.1)
\end{equation}

\begin{figure}[hbtp]
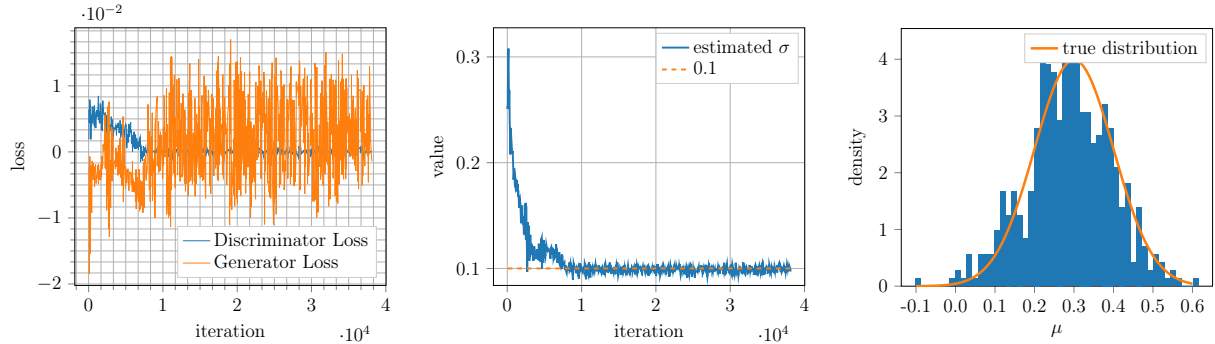

  \scalebox{0.6}{\input{figures/ploss}}~
  \scalebox{0.6}{\input{figures/psigma}}~
  \scalebox{0.6}{\input{figures/poidist}}
  \caption{ANA results for \cref{equ:aaa}. In the first plot, we show the model loss and discriminator loss. In the second plot, we see that the estimated $\sigma$ converges to $\sigma^*$ after around 1000 iterations. The last plot shows the true distribution of $\mu^*$ and generated distribution for $\mu$ after iteration 38,000.}
  \label{fig:poisson}
\end{figure}

\Cref{fig:poisson} shows the results of ANA. In the first plot, we show the model loss and discriminator loss. We see that the discriminator loss converges to 0, implying that the discriminator has been successfully fooled. The model loss oscillates around 0. In the second plot, we see that the estimated $\sigma$ converges to $\sigma^*$ after around 1000 iterations. The values then oscillate around 0.1 due to the intrinsically adversarial optimization. The last plot shows the true distribution of $\mu^*$ and generated distribution for $\mu$ after iteration 38,000. We see that the generated distribution matches the true distribution reasonably well.

\subsection{Uncertainty Quantification with Multimodal Distribution}\label{sect:uqmd}

We consider a variant of \Cref{sect:n1}. This example demonstrates that ANA can be used to learn a rather complicated distribution such as multimodal ones. 

Different from the last section, we assume that $\sigma$ is known. Besides, $\mu$ is subject to a nontrivial distribution, e.g., Gaussian mixture distribution,
\begin{equation}\label{equ:ppp}
    p = 0.4 \, \mathcal{N}(0.3,0.1) + 0.6 \, \mathcal{N}(0.8,0.05)
\end{equation}

\begin{figure}[hbtp]
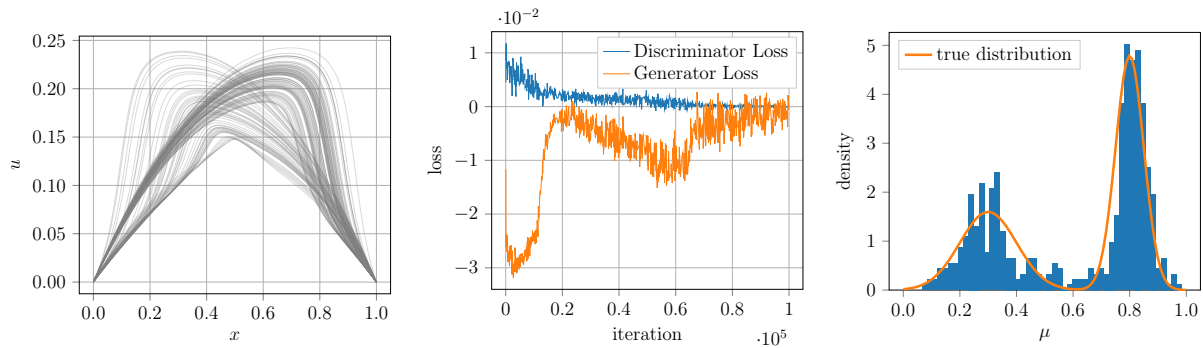

  \scalebox{0.6}{\input{figures/bipoi.tex}}~
  \scalebox{0.6}{\input{figures/bipoiloss.tex}}~
    \scalebox{0.6}{\input{figures/bipoidist.tex}}
  \caption{ANA results for \cref{equ:aaa} with \cref{equ:ppp}. The first plot is generated samples of $u(x)$. The second plot shows the loss functions. The last plot shows the generated distribution at iteration 100,000 together with the true distribution.}
  \label{fig:bipoi}
\end{figure}

The generated samples are shown in the first plot in \cref{fig:bipoi}. ANA is carried out with the same setting as in \cref{sect:n1}. We show the loss functions and the generated distribution in \cref{fig:bipoi}. We can see that the loss functions for $F_\theta$ and the discrimination converge to 0. The generated distribution also captures the multimodal of the distribution correctly. 

\begin{figure}[hbtp]
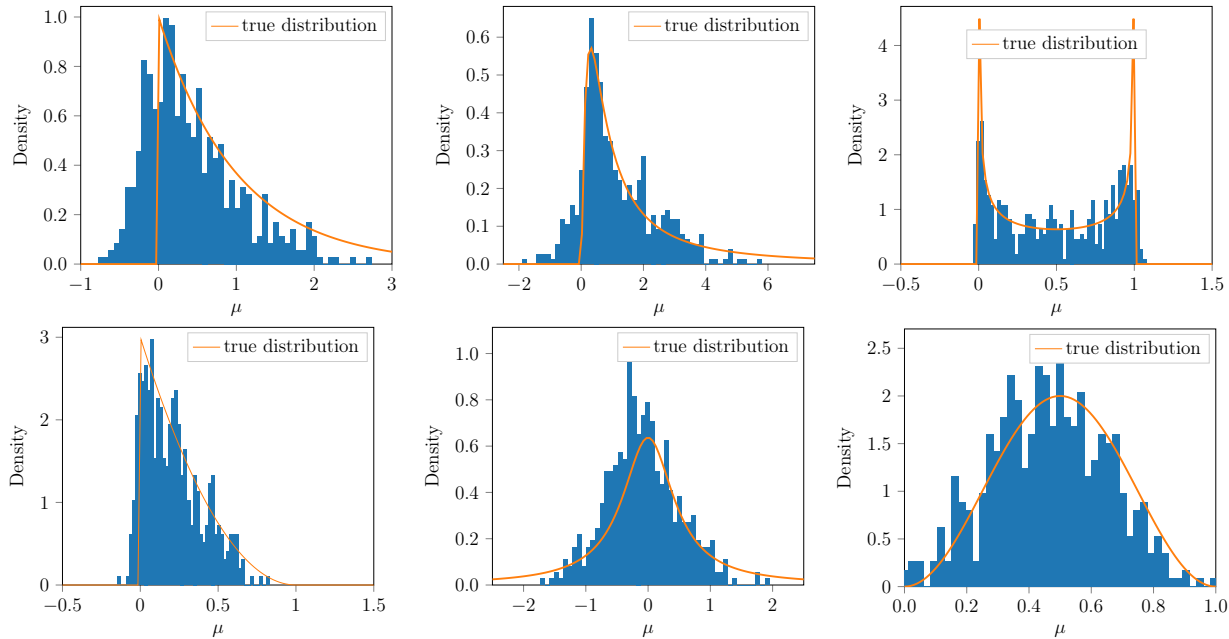

  \scalebox{0.6}{\input{figures/1D/d2_result1}}~
  \scalebox{0.6}{\input{figures/1D/d2_result2}}~
    \scalebox{0.6}{
\begin{tikzpicture}

\definecolor{color1}{rgb}{1,0.498039215686275,0.0549019607843137}
\definecolor{color0}{rgb}{0.12156862745098,0.466666666666667,0.705882352941177}

\begin{axis}[
legend cell align={left},
legend entries={{true distribution}},
legend style={at={(0.5,0.91)}, anchor=north, draw=white!80.0!black},
tick align=outside,
tick pos=left,
x grid style={white!69.01960784313725!black},
xlabel={$\mu$},
xmin=-0.5, xmax=1.5,
y grid style={white!69.01960784313725!black},
ylabel={Density},
ymin=0, ymax=4.71488933345991
]
\addlegendimage{no markers, color1}
\draw[fill=color0,draw opacity=0] (axis cs:-0.0349011449832977,0) rectangle (axis cs:-0.0126641557379147,0.719521866177188);
\draw[fill=color0,draw opacity=0] (axis cs:-0.0126641557379147,0) rectangle (axis cs:0.00957283350746826,2.24850583180371);
\draw[fill=color0,draw opacity=0] (axis cs:0.00957283350746826,0) rectangle (axis cs:0.0318098227528512,2.60826676489231);
\draw[fill=color0,draw opacity=0] (axis cs:0.0318098227528512,0) rectangle (axis cs:0.0540468119982342,1.52898396562652);
\draw[fill=color0,draw opacity=0] (axis cs:0.0540468119982342,0) rectangle (axis cs:0.0762838012436171,1.25916326581008);
\draw[fill=color0,draw opacity=0] (axis cs:0.0762838012436171,0) rectangle (axis cs:0.0985207904890001,1.07928279926578);
\draw[fill=color0,draw opacity=0] (axis cs:0.0985207904890001,0) rectangle (axis cs:0.120757779734383,0.449701166360742);
\draw[fill=color0,draw opacity=0] (axis cs:0.120757779734383,0) rectangle (axis cs:0.142994768979766,1.16922303253793);
\draw[fill=color0,draw opacity=0] (axis cs:0.142994768979766,0) rectangle (axis cs:0.165231758225149,1.07928279926578);
\draw[fill=color0,draw opacity=0] (axis cs:0.165231758225149,0) rectangle (axis cs:0.187468747470532,1.07928279926578);
\draw[fill=color0,draw opacity=0] (axis cs:0.187468747470532,0) rectangle (axis cs:0.209705736715915,0.809462099449337);
\draw[fill=color0,draw opacity=0] (axis cs:0.209705736715915,0) rectangle (axis cs:0.231942725961298,0.539641399632891);
\draw[fill=color0,draw opacity=0] (axis cs:0.231942725961298,0) rectangle (axis cs:0.254179715206681,0.179880466544297);
\draw[fill=color0,draw opacity=0] (axis cs:0.254179715206681,0) rectangle (axis cs:0.276416704452064,0.53964139963289);
\draw[fill=color0,draw opacity=0] (axis cs:0.276416704452064,0) rectangle (axis cs:0.298653693697447,0.539641399632892);
\draw[fill=color0,draw opacity=0] (axis cs:0.298653693697447,0) rectangle (axis cs:0.32089068294283,0.899402332721484);
\draw[fill=color0,draw opacity=0] (axis cs:0.32089068294283,0) rectangle (axis cs:0.343127672188213,0.899402332721484);
\draw[fill=color0,draw opacity=0] (axis cs:0.343127672188213,0) rectangle (axis cs:0.365364661433596,0.809462099449338);
\draw[fill=color0,draw opacity=0] (axis cs:0.365364661433596,0) rectangle (axis cs:0.387601650678979,0.449701166360742);
\draw[fill=color0,draw opacity=0] (axis cs:0.387601650678979,0) rectangle (axis cs:0.409838639924362,0.629581632905039);
\draw[fill=color0,draw opacity=0] (axis cs:0.409838639924362,0) rectangle (axis cs:0.432075629169744,0.539641399632892);
\draw[fill=color0,draw opacity=0] (axis cs:0.432075629169744,0) rectangle (axis cs:0.454312618415127,0.809462099449336);
\draw[fill=color0,draw opacity=0] (axis cs:0.454312618415127,0) rectangle (axis cs:0.47654960766051,1.07928279926578);
\draw[fill=color0,draw opacity=0] (axis cs:0.47654960766051,0) rectangle (axis cs:0.498786596905893,0.809462099449336);
\draw[fill=color0,draw opacity=0] (axis cs:0.498786596905893,0) rectangle (axis cs:0.521023586151276,0.629581632905039);
\draw[fill=color0,draw opacity=0] (axis cs:0.521023586151276,0) rectangle (axis cs:0.543260575396659,0.719521866177187);
\draw[fill=color0,draw opacity=0] (axis cs:0.543260575396659,0) rectangle (axis cs:0.565497564642042,0.0899402332721484);
\draw[fill=color0,draw opacity=0] (axis cs:0.565497564642042,0) rectangle (axis cs:0.587734553887425,0.629581632905039);
\draw[fill=color0,draw opacity=0] (axis cs:0.587734553887425,0) rectangle (axis cs:0.609971543132808,0.989342565993637);
\draw[fill=color0,draw opacity=0] (axis cs:0.609971543132808,0) rectangle (axis cs:0.632208532378191,0.449701166360742);
\draw[fill=color0,draw opacity=0] (axis cs:0.632208532378191,0) rectangle (axis cs:0.654445521623574,0.449701166360742);
\draw[fill=color0,draw opacity=0] (axis cs:0.654445521623574,0) rectangle (axis cs:0.676682510868957,0.53964139963289);
\draw[fill=color0,draw opacity=0] (axis cs:0.676682510868957,0) rectangle (axis cs:0.69891950011434,0.449701166360742);
\draw[fill=color0,draw opacity=0] (axis cs:0.69891950011434,0) rectangle (axis cs:0.721156489359723,0.809462099449336);
\draw[fill=color0,draw opacity=0] (axis cs:0.721156489359723,0) rectangle (axis cs:0.743393478605106,1.16922303253793);
\draw[fill=color0,draw opacity=0] (axis cs:0.743393478605106,0) rectangle (axis cs:0.765630467850489,0.359760933088595);
\draw[fill=color0,draw opacity=0] (axis cs:0.765630467850489,0) rectangle (axis cs:0.787867457095872,0.719521866177187);
\draw[fill=color0,draw opacity=0] (axis cs:0.787867457095872,0) rectangle (axis cs:0.810104446341255,0.269820699816445);
\draw[fill=color0,draw opacity=0] (axis cs:0.810104446341255,0) rectangle (axis cs:0.832341435586638,1.16922303253793);
\draw[fill=color0,draw opacity=0] (axis cs:0.832341435586638,0) rectangle (axis cs:0.854578424832021,0.899402332721484);
\draw[fill=color0,draw opacity=0] (axis cs:0.854578424832021,0) rectangle (axis cs:0.876815414077404,1.43904373235437);
\draw[fill=color0,draw opacity=0] (axis cs:0.876815414077404,0) rectangle (axis cs:0.899052403322787,0.80946209944934);
\draw[fill=color0,draw opacity=0] (axis cs:0.899052403322787,0) rectangle (axis cs:0.92128939256817,1.70886443217082);
\draw[fill=color0,draw opacity=0] (axis cs:0.92128939256817,0) rectangle (axis cs:0.943526381813553,1.79880466544297);
\draw[fill=color0,draw opacity=0] (axis cs:0.943526381813553,0) rectangle (axis cs:0.965763371058936,1.43904373235437);
\draw[fill=color0,draw opacity=0] (axis cs:0.965763371058936,0) rectangle (axis cs:0.988000360304318,1.79880466544298);
\draw[fill=color0,draw opacity=0] (axis cs:0.988000360304318,0) rectangle (axis cs:1.0102373495497,1.16922303253792);
\draw[fill=color0,draw opacity=0] (axis cs:1.0102373495497,0) rectangle (axis cs:1.03247433879508,1.34910349908223);
\draw[fill=color0,draw opacity=0] (axis cs:1.03247433879508,0) rectangle (axis cs:1.05471132804047,0.269820699816444);
\draw[fill=color0,draw opacity=0] (axis cs:1.05471132804047,0) rectangle (axis cs:1.07694831728585,0.0899402332721488);
\addplot [very thick, color1]
table [row sep=\\]{%
-0.5	0 \\
-0.47979797979798	0 \\
-0.45959595959596	0 \\
-0.439393939393939	0 \\
-0.419191919191919	0 \\
-0.398989898989899	0 \\
-0.378787878787879	0 \\
-0.358585858585859	0 \\
-0.338383838383838	0 \\
-0.318181818181818	0 \\
-0.297979797979798	0 \\
-0.277777777777778	0 \\
-0.257575757575758	0 \\
-0.237373737373737	0 \\
-0.217171717171717	0 \\
-0.196969696969697	0 \\
-0.176767676767677	0 \\
-0.156565656565657	0 \\
-0.136363636363636	0 \\
-0.116161616161616	0 \\
-0.095959595959596	0 \\
-0.0757575757575757	0 \\
-0.0555555555555556	0 \\
-0.0353535353535354	0 \\
-0.0151515151515151	0 \\
0.00505050505050508	4.49037079377129 \\
0.0252525252525253	2.02885804259676 \\
0.0454545454545454	1.52814005969429 \\
0.0656565656565656	1.2851617084421 \\
0.0858585858585859	1.13619158765545 \\
0.106060606060606	1.03375836549109 \\
0.126262626262626	0.958345839263334 \\
0.146464646464647	0.900270389852276 \\
0.166666666666667	0.854115052100613 \\
0.186868686868687	0.81658522467735 \\
0.207070707070707	0.785549511708052 \\
0.227272727272727	0.759562202718595 \\
0.247474747474748	0.737605494813596 \\
0.267676767676768	0.718941215611386 \\
0.287878787878788	0.703021057211265 \\
0.308080808080808	0.689429903196885 \\
0.328282828282828	0.677848780035936 \\
0.348484848484849	0.668029924272091 \\
0.368686868686869	0.659779596533234 \\
0.388888888888889	0.652946006201518 \\
0.409090909090909	0.64741070653079 \\
0.429292929292929	0.643082413439515 \\
0.44949494949495	0.639892566773581 \\
0.46969696969697	0.637792185464555 \\
0.48989898989899	0.636749721291647 \\
0.51010101010101	0.636749721291647 \\
0.53030303030303	0.637792185464555 \\
0.550505050505051	0.639892566773581 \\
0.570707070707071	0.643082413439515 \\
0.590909090909091	0.64741070653079 \\
0.611111111111111	0.652946006201518 \\
0.631313131313131	0.659779596533234 \\
0.651515151515152	0.668029924272091 \\
0.671717171717172	0.677848780035936 \\
0.691919191919192	0.689429903196885 \\
0.712121212121212	0.703021057211265 \\
0.732323232323232	0.718941215611386 \\
0.752525252525253	0.737605494813596 \\
0.772727272727273	0.759562202718595 \\
0.792929292929293	0.785549511708052 \\
0.813131313131313	0.81658522467735 \\
0.833333333333333	0.854115052100612 \\
0.853535353535354	0.900270389852276 \\
0.873737373737374	0.958345839263334 \\
0.893939393939394	1.03375836549109 \\
0.914141414141414	1.13619158765545 \\
0.934343434343434	1.2851617084421 \\
0.954545454545454	1.52814005969428 \\
0.974747474747475	2.02885804259676 \\
0.994949494949495	4.49037079377134 \\
1.01515151515152	0 \\
1.03535353535354	0 \\
1.05555555555556	0 \\
1.07575757575758	0 \\
1.0959595959596	0 \\
1.11616161616162	0 \\
1.13636363636364	0 \\
1.15656565656566	0 \\
1.17676767676768	0 \\
1.1969696969697	0 \\
1.21717171717172	0 \\
1.23737373737374	0 \\
1.25757575757576	0 \\
1.27777777777778	0 \\
1.2979797979798	0 \\
1.31818181818182	0 \\
1.33838383838384	0 \\
1.35858585858586	0 \\
1.37878787878788	0 \\
1.3989898989899	0 \\
1.41919191919192	0 \\
1.43939393939394	0 \\
1.45959595959596	0 \\
1.47979797979798	0 \\
1.5	0 \\
};

\end{axis}

\end{tikzpicture}}
    \scalebox{0.6}{
\begin{tikzpicture}

\definecolor{color1}{rgb}{1,0.498039215686275,0.0549019607843137}
\definecolor{color0}{rgb}{0.12156862745098,0.466666666666667,0.705882352941177}

\begin{axis}[
legend cell align={left},
legend entries={{true distribution}},
legend style={draw=white!80.0!black},
tick align=outside,
tick pos=left,
x grid style={white!69.01960784313725!black},
xlabel={$\mu$},
xmin=-0.5, xmax=1.5,
y grid style={white!69.01960784313725!black},
ylabel={Density},
ymin=0, ymax=3.11974234039659
]
\addlegendimage{no markers, color1}
\draw[fill=color0,draw opacity=0] (axis cs:-0.148969736042318,0) rectangle (axis cs:-0.129448893179295,0.102454592459658);
\draw[fill=color0,draw opacity=0] (axis cs:-0.129448893179295,0) rectangle (axis cs:-0.109928050316272,0);
\draw[fill=color0,draw opacity=0] (axis cs:-0.109928050316273,0) rectangle (axis cs:-0.0904072074532499,0);
\draw[fill=color0,draw opacity=0] (axis cs:-0.0904072074532499,0) rectangle (axis cs:-0.0708863645902274,0.102454592459658);
\draw[fill=color0,draw opacity=0] (axis cs:-0.0708863645902274,0) rectangle (axis cs:-0.0513655217272048,0.614727554757948);
\draw[fill=color0,draw opacity=0] (axis cs:-0.0513655217272048,0) rectangle (axis cs:-0.0318446788641823,1.02454592459658);
\draw[fill=color0,draw opacity=0] (axis cs:-0.0318446788641823,0) rectangle (axis cs:-0.0123238360011597,2.04909184919316);
\draw[fill=color0,draw opacity=0] (axis cs:-0.0123238360011597,0) rectangle (axis cs:0.00719700686186286,2.56136481149145);
\draw[fill=color0,draw opacity=0] (axis cs:0.00719700686186286,0) rectangle (axis cs:0.0267178497248854,2.45891021903179);
\draw[fill=color0,draw opacity=0] (axis cs:0.0267178497248854,0) rectangle (axis cs:0.046238692587908,2.66381940395111);
\draw[fill=color0,draw opacity=0] (axis cs:0.046238692587908,0) rectangle (axis cs:0.0657595354509306,2.35645562657213);
\draw[fill=color0,draw opacity=0] (axis cs:0.0657595354509306,0) rectangle (axis cs:0.0852803783139531,2.97118318133008);
\draw[fill=color0,draw opacity=0] (axis cs:0.0852803783139531,0) rectangle (axis cs:0.104801221176976,1.53681888689487);
\draw[fill=color0,draw opacity=0] (axis cs:0.104801221176976,0) rectangle (axis cs:0.124322064039998,2.25400103411247);
\draw[fill=color0,draw opacity=0] (axis cs:0.124322064039998,0) rectangle (axis cs:0.143842906903021,2.15154644165282);
\draw[fill=color0,draw opacity=0] (axis cs:0.143842906903021,0) rectangle (axis cs:0.163363749766043,1.53681888689487);
\draw[fill=color0,draw opacity=0] (axis cs:0.163363749766043,0) rectangle (axis cs:0.182884592629066,1.43436429443521);
\draw[fill=color0,draw opacity=0] (axis cs:0.182884592629066,0) rectangle (axis cs:0.202405435492088,1.94663725673351);
\draw[fill=color0,draw opacity=0] (axis cs:0.202405435492088,0) rectangle (axis cs:0.221926278355111,2.25400103411247);
\draw[fill=color0,draw opacity=0] (axis cs:0.221926278355111,0) rectangle (axis cs:0.241447121218134,2.35645562657214);
\draw[fill=color0,draw opacity=0] (axis cs:0.241447121218134,0) rectangle (axis cs:0.260967964081156,1.9466372567335);
\draw[fill=color0,draw opacity=0] (axis cs:0.260967964081156,0) rectangle (axis cs:0.280488806944179,1.33190970197555);
\draw[fill=color0,draw opacity=0] (axis cs:0.280488806944179,0) rectangle (axis cs:0.300009649807201,1.63927347935453);
\draw[fill=color0,draw opacity=0] (axis cs:0.300009649807201,0) rectangle (axis cs:0.319530492670224,1.02454592459658);
\draw[fill=color0,draw opacity=0] (axis cs:0.319530492670224,0) rectangle (axis cs:0.339051335533246,0.717182147217605);
\draw[fill=color0,draw opacity=0] (axis cs:0.339051335533246,0) rectangle (axis cs:0.358572178396269,1.33190970197556);
\draw[fill=color0,draw opacity=0] (axis cs:0.358572178396269,0) rectangle (axis cs:0.378093021259292,1.12700051705624);
\draw[fill=color0,draw opacity=0] (axis cs:0.378093021259292,0) rectangle (axis cs:0.397613864122314,0.614727554757947);
\draw[fill=color0,draw opacity=0] (axis cs:0.397613864122314,0) rectangle (axis cs:0.417134706985337,0.512272962298289);
\draw[fill=color0,draw opacity=0] (axis cs:0.417134706985337,0) rectangle (axis cs:0.436655549848359,0.717182147217609);
\draw[fill=color0,draw opacity=0] (axis cs:0.436655549848359,0) rectangle (axis cs:0.456176392711382,1.22945510951589);
\draw[fill=color0,draw opacity=0] (axis cs:0.456176392711382,0) rectangle (axis cs:0.475697235574404,1.33190970197555);
\draw[fill=color0,draw opacity=0] (axis cs:0.475697235574404,0) rectangle (axis cs:0.495218078437427,0.614727554757947);
\draw[fill=color0,draw opacity=0] (axis cs:0.495218078437427,0) rectangle (axis cs:0.514738921300449,0.717182147217607);
\draw[fill=color0,draw opacity=0] (axis cs:0.514738921300449,0) rectangle (axis cs:0.534259764163472,0.614727554757951);
\draw[fill=color0,draw opacity=0] (axis cs:0.534259764163472,0) rectangle (axis cs:0.553780607026495,0.307363777378972);
\draw[fill=color0,draw opacity=0] (axis cs:0.553780607026495,0) rectangle (axis cs:0.573301449889517,0.409818369838634);
\draw[fill=color0,draw opacity=0] (axis cs:0.573301449889517,0) rectangle (axis cs:0.59282229275254,0.614727554757944);
\draw[fill=color0,draw opacity=0] (axis cs:0.59282229275254,0) rectangle (axis cs:0.612343135615562,0.614727554757951);
\draw[fill=color0,draw opacity=0] (axis cs:0.612343135615562,0) rectangle (axis cs:0.631863978478585,0.614727554757951);
\draw[fill=color0,draw opacity=0] (axis cs:0.631863978478585,0) rectangle (axis cs:0.651384821341607,0.204909184919315);
\draw[fill=color0,draw opacity=0] (axis cs:0.651384821341607,0) rectangle (axis cs:0.67090566420463,0.307363777378975);
\draw[fill=color0,draw opacity=0] (axis cs:0.67090566420463,0) rectangle (axis cs:0.690426507067653,0);
\draw[fill=color0,draw opacity=0] (axis cs:0.690426507067653,0) rectangle (axis cs:0.709947349930675,0);
\draw[fill=color0,draw opacity=0] (axis cs:0.709947349930675,0) rectangle (axis cs:0.729468192793698,0.102454592459657);
\draw[fill=color0,draw opacity=0] (axis cs:0.729468192793698,0) rectangle (axis cs:0.74898903565672,0);
\draw[fill=color0,draw opacity=0] (axis cs:0.74898903565672,0) rectangle (axis cs:0.768509878519743,0);
\draw[fill=color0,draw opacity=0] (axis cs:0.768509878519743,0) rectangle (axis cs:0.788030721382765,0.102454592459657);
\draw[fill=color0,draw opacity=0] (axis cs:0.788030721382765,0) rectangle (axis cs:0.807551564245788,0);
\draw[fill=color0,draw opacity=0] (axis cs:0.807551564245788,0) rectangle (axis cs:0.82707240710881,0.102454592459658);
\addplot [semithick, color1]
table [row sep=\\]{%
-0.5	0 \\
-0.47979797979798	0 \\
-0.45959595959596	0 \\
-0.439393939393939	0 \\
-0.419191919191919	0 \\
-0.398989898989899	0 \\
-0.378787878787879	0 \\
-0.358585858585859	0 \\
-0.338383838383838	0 \\
-0.318181818181818	0 \\
-0.297979797979798	0 \\
-0.277777777777778	0 \\
-0.257575757575758	0 \\
-0.237373737373737	0 \\
-0.217171717171717	0 \\
-0.196969696969697	0 \\
-0.176767676767677	0 \\
-0.156565656565657	0 \\
-0.136363636363636	0 \\
-0.116161616161616	0 \\
-0.095959595959596	0 \\
-0.0757575757575757	0 \\
-0.0555555555555556	0 \\
-0.0353535353535354	0 \\
-0.0151515151515151	0 \\
0.00505050505050508	2.96977349250077 \\
0.0252525252525253	2.85039791857974 \\
0.0454545454545454	2.73347107438017 \\
0.0656565656565656	2.61899295990205 \\
0.0858585858585859	2.50696357514539 \\
0.106060606060606	2.39738292011019 \\
0.126262626262626	2.29025099479645 \\
0.146464646464647	2.18556779920416 \\
0.166666666666667	2.08333333333333 \\
0.186868686868687	1.98354759718396 \\
0.207070707070707	1.88621059075605 \\
0.227272727272727	1.79132231404959 \\
0.247474747474748	1.69888276706459 \\
0.267676767676768	1.60889194980104 \\
0.287878787878788	1.52134986225895 \\
0.308080808080808	1.43625650443832 \\
0.328282828282828	1.35361187633915 \\
0.348484848484849	1.27341597796143 \\
0.368686868686869	1.19566880930517 \\
0.388888888888889	1.12037037037037 \\
0.409090909090909	1.04752066115702 \\
0.429292929292929	0.977119681665136 \\
0.44949494949495	0.909167431894705 \\
0.46969696969697	0.84366391184573 \\
0.48989898989899	0.780609121518212 \\
0.51010101010101	0.720003060912152 \\
0.53030303030303	0.661845730027548 \\
0.550505050505051	0.606137128864401 \\
0.570707070707071	0.552877257422712 \\
0.590909090909091	0.50206611570248 \\
0.611111111111111	0.453703703703704 \\
0.631313131313131	0.407790021426385 \\
0.651515151515152	0.364325068870523 \\
0.671717171717172	0.323308846036119 \\
0.691919191919192	0.284741352923171 \\
0.712121212121212	0.24862258953168 \\
0.732323232323232	0.214952555861647 \\
0.752525252525253	0.18373125191307 \\
0.772727272727273	0.15495867768595 \\
0.792929292929293	0.128634833180288 \\
0.813131313131313	0.104759718396082 \\
0.833333333333333	0.0833333333333334 \\
0.853535353535354	0.0643556779920416 \\
0.873737373737374	0.047826752372207 \\
0.893939393939394	0.0337465564738292 \\
0.914141414141414	0.0221150902969085 \\
0.934343434343434	0.0129323538414448 \\
0.954545454545454	0.00619834710743804 \\
0.974747474747475	0.00191307009488828 \\
0.994949494949495	7.65228037955286e-05 \\
1.01515151515152	0 \\
1.03535353535354	0 \\
1.05555555555556	0 \\
1.07575757575758	0 \\
1.0959595959596	0 \\
1.11616161616162	0 \\
1.13636363636364	0 \\
1.15656565656566	0 \\
1.17676767676768	0 \\
1.1969696969697	0 \\
1.21717171717172	0 \\
1.23737373737374	0 \\
1.25757575757576	0 \\
1.27777777777778	0 \\
1.2979797979798	0 \\
1.31818181818182	0 \\
1.33838383838384	0 \\
1.35858585858586	0 \\
1.37878787878788	0 \\
1.3989898989899	0 \\
1.41919191919192	0 \\
1.43939393939394	0 \\
1.45959595959596	0 \\
1.47979797979798	0 \\
1.5	0 \\
};

\end{axis}

\end{tikzpicture}}~
  \scalebox{0.6}{\input{figures/1D/d2_result5}}~
    \scalebox{0.6}{\input{figures/1D/d2_result7}}
  \caption{ANA results for estimating the probability distribution of $\mu$ where $\mu$ is subject to a different distribution for each case. The first rows are the Exponential distribution, the F distribution, and the Arcsine distribution; the second rows are the Beta distribution, the Cauchy distribution, and the raised Cosine distribution.}
  \label{fig:bipoi2}
\end{figure}

Other distributions for $\mu$ are also tested with the same neural network architecture and optimizer. The results in \Cref{fig:bipoi2} demonstrate that the proposed method is effective and generic in estimating the unknown distributions in the PDE system. The test distributions include
\begin{enumerate}
    \item the Exponential distribution with a rate 1;
    \item the F distribution with degrees of freedom $(5,2)$;
    \item the Arcsine distribution in $[0,1]$;
    \item the Beta distribution with shape parameters $(1,3)$;
    \item the Cauchy distribution with the location parameter $x_0=0$ and the scale parameter $b=0.5$;
    \item the raised Cosine distribution with parameters $\mu=0.5$, $s=0.5$.
\end{enumerate}

\begin{figure}[hbt]
  \includegraphics[width=1.0\textwidth]{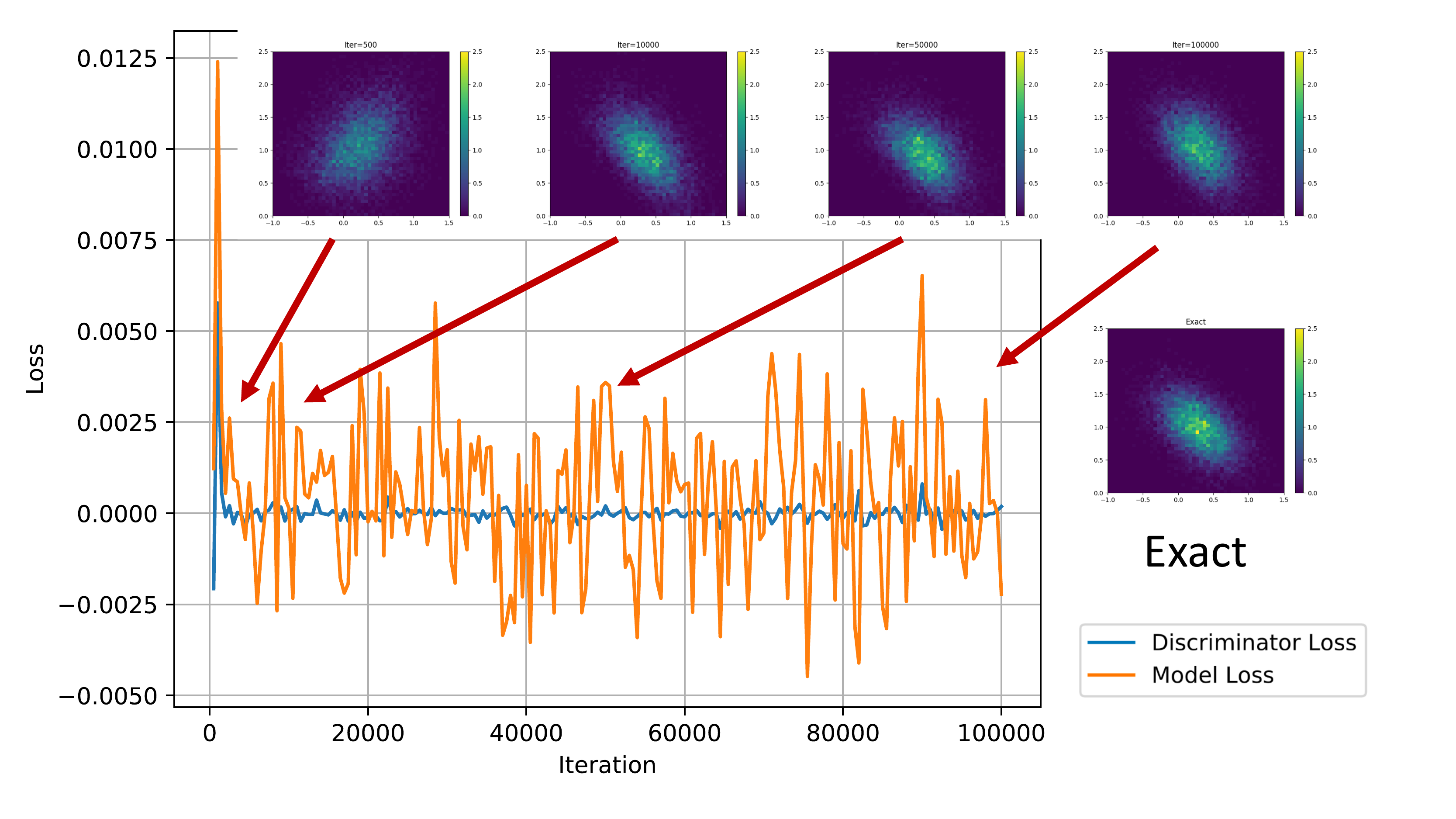}
  \caption{ANA for estimating the distribution of $(\mu, |\sigma|)$ where $(\mu, \sigma)$ is subject to a 2D Gaussian distribution.}
  \label{fig:2D_1}
\end{figure}

\begin{figure}[hbt]
  \includegraphics[width=1.0\textwidth]{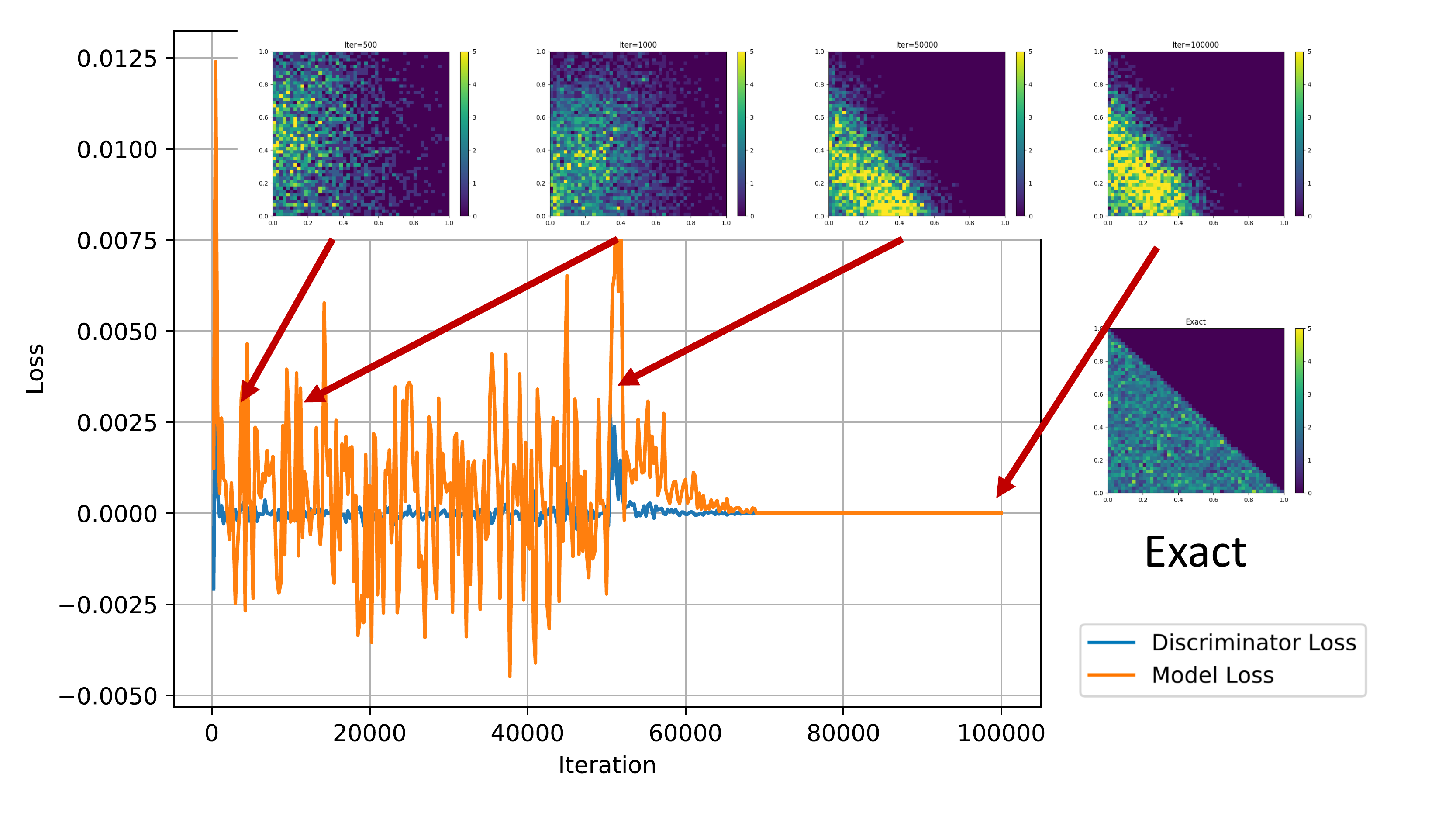}
  \caption{ANA for estimating the distribution of $(\mu, |\sigma|)$ where $(\mu, \sigma)$ is subject to a Dirichlet distribution with parameters $\bm{\alpha}=(1.0, 1.0, 1.0)$.}
  \label{fig:2D_2}
\end{figure}

\begin{figure}[hbt]
  \includegraphics[width=1.0\textwidth]{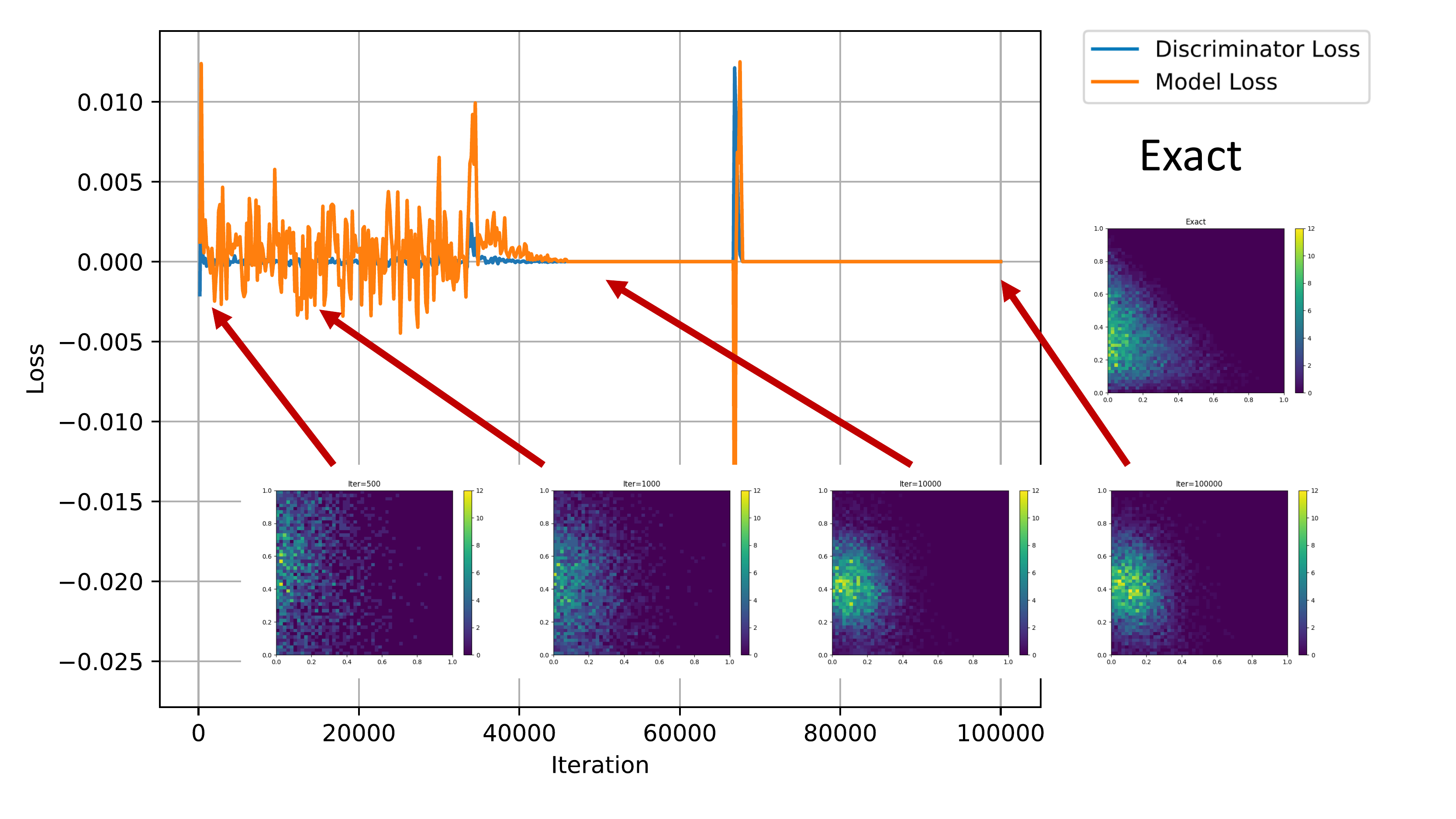}
  \caption{ANA for estimating the distribution of $(\mu, |\sigma|)$ where $(\mu, \sigma)$ is subject to a Dirichlet distribution with parameters $\bm{\alpha}=(1.0,2.0,3.0)$.}
  \label{fig:2D_3}
\end{figure}

We also consider 2D distributions. In this case, we assume that $\mu$, $\sigma$ are both random variables and they are subject to a 2D unknown distribution ($\sigma$ and $\mu$ not necessarily independent). We increase the number of the output layer neurons from 1 to 2 for generating $\sigma$ in the previous case and perform adversarial numerical analysis. For generating synthetic data, we draw $(\mu, \sigma)$ from the following random distributions
\begin{enumerate}
    \item a 2D Gaussian distribution with mean $(0,3,1.0)$ and a covariance matrix $\begin{pmatrix}
        0.1 & -0.05\\
        -0.05 & 0.1
    \end{pmatrix}$;
    \item a Dirichlet distribution with parameters $\bm{\alpha}=(1.0, 1.0, 1.0)$;
    \item a Dirichlet distribution with parameters $\bm{\alpha}=(1.0,2.0,3.0)$.
\end{enumerate} 
since in the model $\sigma$ appears as $\sigma^2$ in the coefficient function $a(x)$, we consider estimating the distribution of $(\mu, |\sigma|)$ instead of $(\mu, \sigma)$. \Cref{fig:2D_1,fig:2D_2,fig:2D_3} show the results at multiple checkpoints. It is remarkable that with the same neural network architecture and optimization scheme, ANA is able to learn a quite similar pattern as the true distribution.


\subsection{Parameter Calibration for CIR Processes from Sample Paths}\label{sect:num_cir}

In this section, we consider the CIR process example considered in \Cref{sect:def,sect:cireg} and numerically demonstrate the analysis in \Cref{sect:ca}. We also compare different optimizers and found that \texttt{LBFGS} converges most smoothly and fastest. 

\begin{figure}[hbtp]
  \scalebox{0.8}{\input{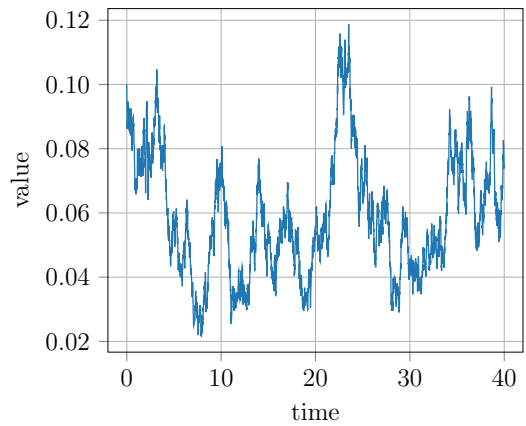}}
  \caption{A sample path of the CIR process.}
  \label{fig:cirpath}
\end{figure}

For better accuracy, we use the weighted Milstein scheme~\cite{glasserman2005large} for the simulation \cref{equ:sim0}
\begin{equation}\label{equ:sim1}
    y = \frac{{x + \kappa (\tau  - \alpha x)\Delta t + \sigma \sqrt x \sqrt {\Delta t} W + \frac{1}{4}{\sigma ^2}\Delta t({W^2} - 1)}}{{1 + (1 - \alpha )\kappa \Delta t}}
\end{equation}
where $\alpha\in [0,1]$ is the weight, $x$ is the sample at last time step. In the numerical example, we let $\alpha=0.5$, $\Delta t= 0.01$, $\sigma=0.08$, $\kappa = 0.5$ and the exact $\tau^*=0.06$. The length of the sample path is $4000$. \Cref{fig:cirpath} shows a realization of the sample path. The KL GAN loss functions and three optimizers~(\texttt{ADAM}, \texttt{RMSProp} and  \texttt{LBFGS}) with full batch size are used. 

\begin{figure}[hbtp]
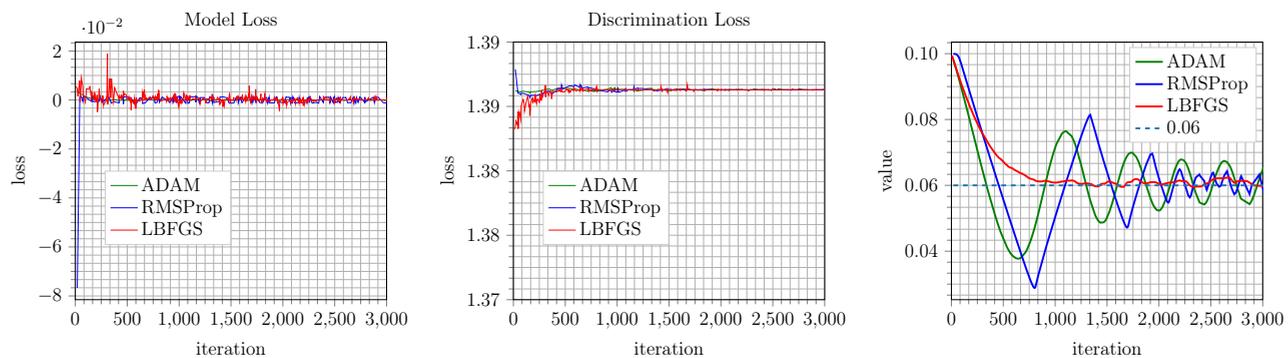

  \scalebox{0.6}{\input{figures/cirg}}~
  \scalebox{0.6}{\input{figures/cird}}~
    \scalebox{0.6}{\input{figures/cirtheta}}
  \caption{ANA results for \cref{equ:sim1}. The first two plots show the model loss and discrimination loss respectively. The third plot shows the convergence of the hidden parameter $\tau$.}
  \label{fig:cirresult}
\end{figure}

\Cref{fig:cirresult} shows the result of ANA. The first two plots show the generator and discrimination losses, and as we can see, they converge to $0$ and $\log 4$ respectively. The third plot shows the convergence profile for $\tau$. It is interesting to see the performance of the different optimizers. With a strong optimizer for the generator~(\texttt{LBFGS}), the convergence profile is much smoother than the others. \texttt{RMSProp} tends to produce kinks and oscillates around the true value $\tau^*=0.06$ like \texttt{ADAM}. 

\begin{figure}[hbtp]
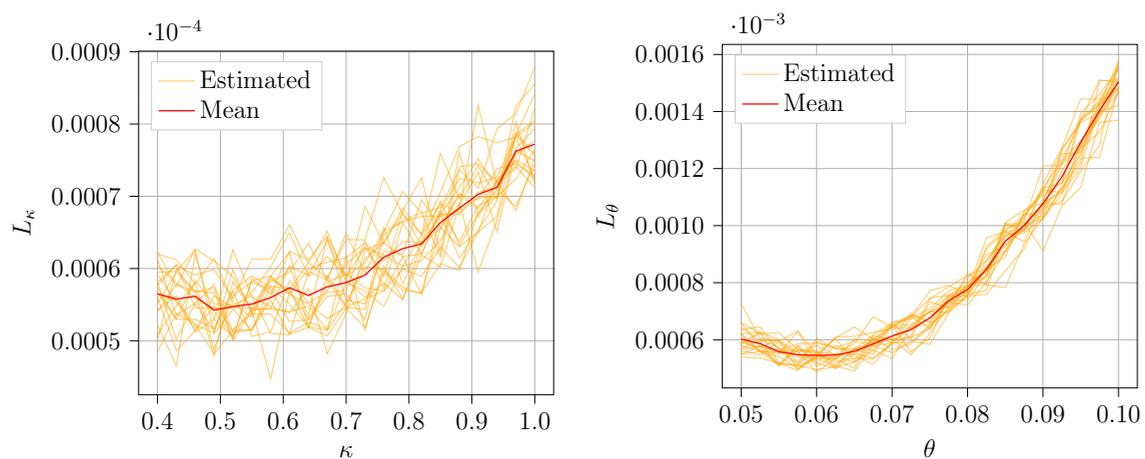

  \scalebox{0.8}{\input{figures/kappa0.tex}}~
  \scalebox{0.8}{\input{figures/theta0.tex}}
  \caption{Discrete KL divergence as a function of $\kappa$ and $\tau$. Each plot has 10 realizations of sample paths. We see that for $\kappa$ the profile is much more oscillating and the landscape at $\kappa^*=0.5$ is nearly flat.}
  \label{fig:kt}
\end{figure}

Next, we discuss the convergence of $\kappa$. As indicated in \Cref{thm:div}, the convergence of  $\kappa$ can be very slow if $X_0=\tau$ and $X_{-1}=\frac{1}{\tau}$, which unfortunately are true if $(R_k$, $R_{k+1})$ are sampled from a long enough sample path. This is demonstrated in \cref{fig:kt}. We approximate the KL divergence with discrete KL divergence 
\begin{equation}
    L_\tau =  \sum_{i} P_{\tau^*, i}\log\left(\frac{P_{\tau^*, i}}{P_{\tau, i}}\right)\qquad L_\kappa =  \sum_{i} P_{\kappa^*, i}\log\left(\frac{P_{\kappa^*, i}}{P_{\kappa, i}}\right)
\end{equation}
where $\tau$, $\kappa$ are the parameters and $\tau^*$, $\kappa^*$ are their true values. $P_{\tau, i}$  is the discrete density function. 
In the first plot, we let $\tau=0.06$, $\sigma=0.08$ and vary $\kappa$. The true value is $\kappa^*=0.5$. In the second plot, we fix $\kappa=0.5$, $\sigma=0.08$ and vary $\tau$, whose true value is $\tau^*=0.06$. We carry out the simulation \cref{equ:sim1} from initial location $R_0=0.05$, $\Delta t = 0.001$ and path length 100,000. The result shows that near the true value, the $L_\kappa$ curve is very flat, indicating small second order derivatives. Compared to $L_\tau$, $L_\kappa$ is much more oscillatory. 

\begin{figure}[hbtp]
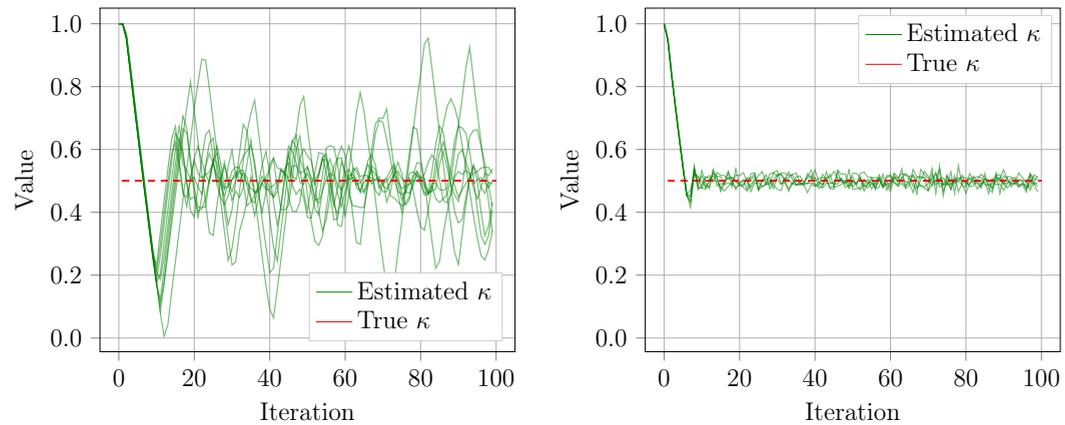

  \scalebox{0.8}{\input{figures/loss2.tex}}~
  \scalebox{0.8}{\input{figures/loss1.tex}}
  \caption{Convergence of $\kappa$ to true value $\kappa^*$.}
  \label{fig:loss12}
\end{figure}

A possible fix to the problem is by distorting the distribution $R_k$. In \cref{fig:loss12}-left, we artificially sample $R_k$ from $\mathcal{U}(0.001, 0.03)$. We use $\Delta t = 0.001$, \texttt{RMSProp} with learning rate $10^{-3}$ for both $\kappa$ and $D_{\xi}$. In each step, $D_{\xi}$ is updated five times while $F_\theta$ is updated only once. We can see that the estimated $\kappa$ oscillates around $\kappa^*=0.5$ but the variance is still large. In \cref{fig:loss12}-right, we further consider reducing the variance by increasing $n$, i.e., the number of samples. Equivalently, we sample $(R_k$, $R_{k+1})$ with time interval $\Delta t=0.1$. 

\subsection{Volatility Inference from Option Prices: Direct Estimation from Monte Carlo Simulation}

The final example is an application of ANA to estimation volatility from European call option prices. In this case, $\sigma$ is the unknown volatility parameter. 

\begin{figure}[hbtp]
  \includegraphics[width=0.8\textwidth]{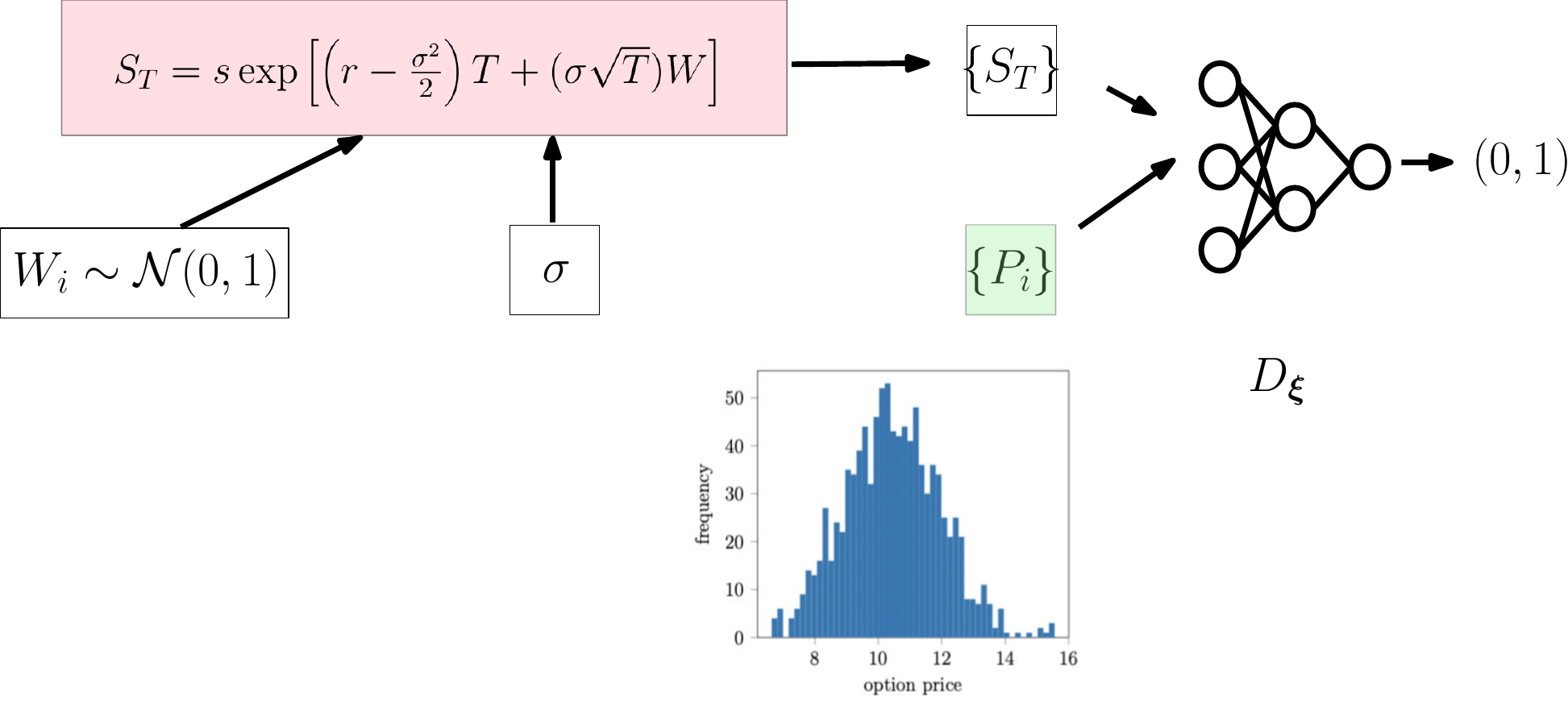}
  \caption{ANA computation pipeline for the option price example.}
  \label{fig:bspip}
\end{figure}

A European call option is a contract giving the holder the right to buy the underlying asset for a fixed strike price $K$ at expiry time $T$. Let $S_t$ be the price of the underlying such as the stock price at time $t$, then the payoff of a European call is a random variable~\cite{luenberger1997investment}
\begin{equation}
P =    \begin{cases}
        S_T - K & \mbox{ if } S_T>K\\
        0 & \mbox{ otherwise }
    \end{cases}
\end{equation}

The stock price can bemodeled as a stochastic process such as geometric Brownian motion
\begin{equation}
    dS_t = r S_tdt + \sigma S_t dW_t\qquad S_0 = s
\end{equation}
where $W_t$ is a Gaussian random variable with zero mean, $\sigma$ is the volatility, $\mu$ is the expected return, and $s$ is the spot price of the stock. 

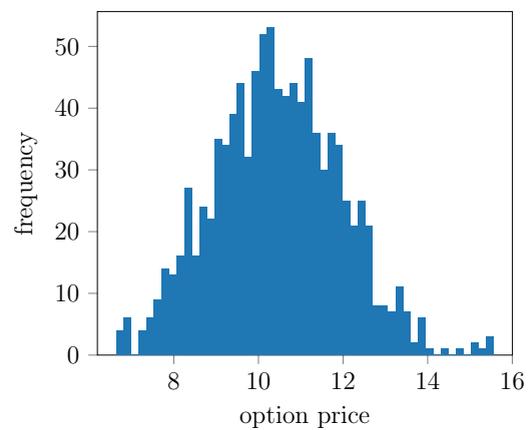
\begin{figure}[hbtp]
  \scalebox{0.8}{
\begin{tikzpicture}

\definecolor{color0}{rgb}{0.12156862745098,0.466666666666667,0.705882352941177}

\begin{axis}[
tick align=outside,
tick pos=left,
x grid style={white!69.01960784313725!black},
xlabel={option price},
xmin=6.19279994036245, xmax=16.0038181588964,
y grid style={white!69.01960784313725!black},
ylabel={frequency},
ymin=0, ymax=55.65
]
\draw[fill=color0,draw opacity=0] (axis cs:6.63875531393218,0) rectangle (axis cs:6.81713746336007,4);
\draw[fill=color0,draw opacity=0] (axis cs:6.81713746336007,0) rectangle (axis cs:6.99551961278796,6);
\draw[fill=color0,draw opacity=0] (axis cs:6.99551961278796,0) rectangle (axis cs:7.17390176221585,0);
\draw[fill=color0,draw opacity=0] (axis cs:7.17390176221585,0) rectangle (axis cs:7.35228391164374,4);
\draw[fill=color0,draw opacity=0] (axis cs:7.35228391164374,0) rectangle (axis cs:7.53066606107163,6);
\draw[fill=color0,draw opacity=0] (axis cs:7.53066606107163,0) rectangle (axis cs:7.70904821049952,9);
\draw[fill=color0,draw opacity=0] (axis cs:7.70904821049952,0) rectangle (axis cs:7.88743035992741,14);
\draw[fill=color0,draw opacity=0] (axis cs:7.88743035992741,0) rectangle (axis cs:8.0658125093553,13);
\draw[fill=color0,draw opacity=0] (axis cs:8.0658125093553,0) rectangle (axis cs:8.24419465878319,16);
\draw[fill=color0,draw opacity=0] (axis cs:8.24419465878319,0) rectangle (axis cs:8.42257680821108,27);
\draw[fill=color0,draw opacity=0] (axis cs:8.42257680821108,0) rectangle (axis cs:8.60095895763897,16);
\draw[fill=color0,draw opacity=0] (axis cs:8.60095895763897,0) rectangle (axis cs:8.77934110706686,24);
\draw[fill=color0,draw opacity=0] (axis cs:8.77934110706686,0) rectangle (axis cs:8.95772325649475,22);
\draw[fill=color0,draw opacity=0] (axis cs:8.95772325649475,0) rectangle (axis cs:9.13610540592264,35);
\draw[fill=color0,draw opacity=0] (axis cs:9.13610540592264,0) rectangle (axis cs:9.31448755535053,34);
\draw[fill=color0,draw opacity=0] (axis cs:9.31448755535053,0) rectangle (axis cs:9.49286970477842,39);
\draw[fill=color0,draw opacity=0] (axis cs:9.49286970477842,0) rectangle (axis cs:9.67125185420631,44);
\draw[fill=color0,draw opacity=0] (axis cs:9.67125185420631,0) rectangle (axis cs:9.8496340036342,32);
\draw[fill=color0,draw opacity=0] (axis cs:9.8496340036342,0) rectangle (axis cs:10.0280161530621,46);
\draw[fill=color0,draw opacity=0] (axis cs:10.0280161530621,0) rectangle (axis cs:10.20639830249,52);
\draw[fill=color0,draw opacity=0] (axis cs:10.20639830249,0) rectangle (axis cs:10.3847804519179,53);
\draw[fill=color0,draw opacity=0] (axis cs:10.3847804519179,0) rectangle (axis cs:10.5631626013458,43);
\draw[fill=color0,draw opacity=0] (axis cs:10.5631626013458,0) rectangle (axis cs:10.7415447507736,42);
\draw[fill=color0,draw opacity=0] (axis cs:10.7415447507736,0) rectangle (axis cs:10.9199269002015,44);
\draw[fill=color0,draw opacity=0] (axis cs:10.9199269002015,0) rectangle (axis cs:11.0983090496294,41);
\draw[fill=color0,draw opacity=0] (axis cs:11.0983090496294,0) rectangle (axis cs:11.2766911990573,48);
\draw[fill=color0,draw opacity=0] (axis cs:11.2766911990573,0) rectangle (axis cs:11.4550733484852,36);
\draw[fill=color0,draw opacity=0] (axis cs:11.4550733484852,0) rectangle (axis cs:11.6334554979131,30);
\draw[fill=color0,draw opacity=0] (axis cs:11.6334554979131,0) rectangle (axis cs:11.811837647341,36);
\draw[fill=color0,draw opacity=0] (axis cs:11.811837647341,0) rectangle (axis cs:11.9902197967689,34);
\draw[fill=color0,draw opacity=0] (axis cs:11.9902197967689,0) rectangle (axis cs:12.1686019461968,25);
\draw[fill=color0,draw opacity=0] (axis cs:12.1686019461968,0) rectangle (axis cs:12.3469840956247,21);
\draw[fill=color0,draw opacity=0] (axis cs:12.3469840956247,0) rectangle (axis cs:12.5253662450525,25);
\draw[fill=color0,draw opacity=0] (axis cs:12.5253662450525,0) rectangle (axis cs:12.7037483944804,21);
\draw[fill=color0,draw opacity=0] (axis cs:12.7037483944804,0) rectangle (axis cs:12.8821305439083,8);
\draw[fill=color0,draw opacity=0] (axis cs:12.8821305439083,0) rectangle (axis cs:13.0605126933362,8);
\draw[fill=color0,draw opacity=0] (axis cs:13.0605126933362,0) rectangle (axis cs:13.2388948427641,7);
\draw[fill=color0,draw opacity=0] (axis cs:13.2388948427641,0) rectangle (axis cs:13.417276992192,11);
\draw[fill=color0,draw opacity=0] (axis cs:13.417276992192,0) rectangle (axis cs:13.5956591416199,7);
\draw[fill=color0,draw opacity=0] (axis cs:13.5956591416199,0) rectangle (axis cs:13.7740412910478,2);
\draw[fill=color0,draw opacity=0] (axis cs:13.7740412910478,0) rectangle (axis cs:13.9524234404757,6);
\draw[fill=color0,draw opacity=0] (axis cs:13.9524234404757,0) rectangle (axis cs:14.1308055899036,1);
\draw[fill=color0,draw opacity=0] (axis cs:14.1308055899036,0) rectangle (axis cs:14.3091877393314,0);
\draw[fill=color0,draw opacity=0] (axis cs:14.3091877393314,0) rectangle (axis cs:14.4875698887593,1);
\draw[fill=color0,draw opacity=0] (axis cs:14.4875698887593,0) rectangle (axis cs:14.6659520381872,0);
\draw[fill=color0,draw opacity=0] (axis cs:14.6659520381872,0) rectangle (axis cs:14.8443341876151,1);
\draw[fill=color0,draw opacity=0] (axis cs:14.8443341876151,0) rectangle (axis cs:15.022716337043,0);
\draw[fill=color0,draw opacity=0] (axis cs:15.022716337043,0) rectangle (axis cs:15.2010984864709,2);
\draw[fill=color0,draw opacity=0] (axis cs:15.2010984864709,0) rectangle (axis cs:15.3794806358988,1);
\draw[fill=color0,draw opacity=0] (axis cs:15.3794806358988,0) rectangle (axis cs:15.5578627853267,3);
\end{axis}

\end{tikzpicture}}
  \caption{Observations of option prices.}
  \label{fig:bshist}
\end{figure}

The task is to estimate $\sigma$ from many observations of option prices $P_i$, $i=1$, $2$, $\ldots$, $100$. We consider the Monte Carlo simulation for the forward problem, 
\begin{equation}
    S_{T} = s\exp\left[ \left(r-\frac{\sigma^2}{2}\right)T+(\sigma\sqrt{T})W \right]
\end{equation}
where $W$ sampled from a standard normal distribution. \Cref{fig:bshist} shows samples for $s=100$, $K=100$, $r=0.05$, $\sigma=0.2$, $T=1$. 

\begin{figure}[hbtp]
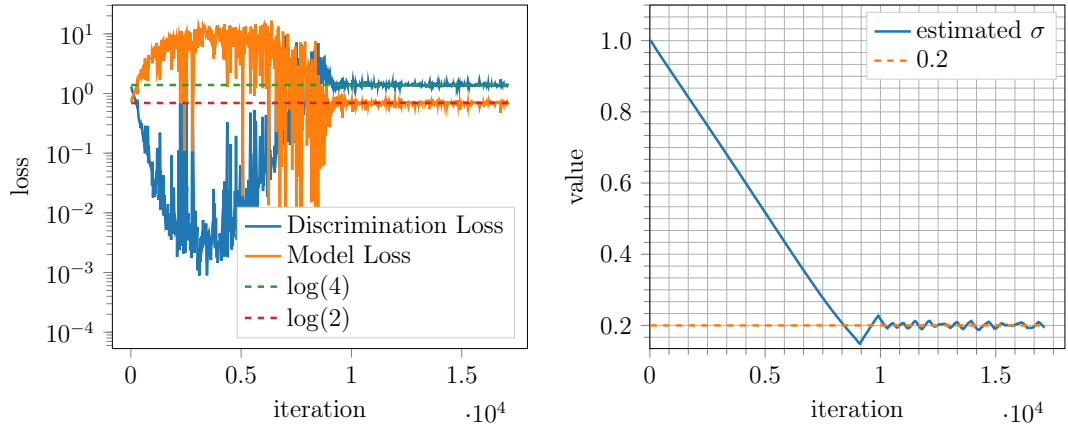

  \scalebox{0.8}{\input{figures/bsloss}}~
  \scalebox{0.8}{\input{figures/bsv}}
  \caption{ANA results for the option pricing example. The first plot shows the loss function while the second shows the convergence of the hidden parameter $\sigma$.}
  \label{fig:bs}
\end{figure}

For ANA, we use the loss functions from the vanilla GAN. \texttt{RMSProp} with learning rate $10^{-4}$ and full batch size is used. The computation pipeline is shown in \cref{fig:bspip}. The result is shown in \cref{fig:bs}. The left plot shows the convergence of the model loss and discrimination loss. We can see that the values converge to theoretical optimal values. The right plot shows the convergence of the estimated volatility.

\section{Conclusion}\label{sect:conc}

We have demonstrated the viability of adversarial numerical analysis for solving inverse modeling problems of the type
\begin{equation}
    x = F(w, \theta)
\end{equation}
where $w$ is a known stochastic process, $\theta$ is an unknown random variable and the output $x$ is also a random variable or process. The key ingredients for applying ANA is: 
\begin{itemize}
    \item Approximating the unknown distribution of $\theta$ by a neural network, $\tilde \theta = G_{\eta}(\theta)$.
    \item Using a discriminator neural network $D_{\xi}$ to measure the discrepancy between the actual random process $x$ and estimated random process $\tilde x = F(w, \tilde \theta)$.
    \item Choosing a proper loss function for $L^D$ and $L^F$. 
    \item Adversarially training $\tilde\theta$ and $\xi$. The gradients $\frac{\partial L^D}{\partial \xi}$ and $\frac{\partial L^F}{\partial \eta}$ can be computed with automatic differentiation. 
\end{itemize}

We also present a computing framework, \texttt{ADCME}, for implementing ANA. We demonstrate numerically that the approach is able to learn underlying parameters and recover complex unknown distributions. \ed{you should summarize your main benchmarks}  Different optimizers were benchmarked and the results showed the promising behavior of full-batch optimizers such a \texttt{LBFGS} for engineering applications. 

However, the approach also has some limitations. For example, the training of GANs can be quite expensive. To obtain the results in \cref{fig:bipoi}, we have run 100,000 iterations, which takes up several CPU hours. Finding more efficient training algorithms, designing specialized loss functions and implementing easy-to-use high-performance toolsets is a promising line of research in the future. 

\newpage

\appendix
\section{Probability Metrics for ANA}\label{sect:pm}

Here we discuss a few probability metrics that are related to ANA; we refer the readers to \cite{gibbs2002choosing} for more details on this topic. Note that these metrics are not mutual exclusive. 

\textit{Information-theoretical Metrics}~\cite{crooks2017measures} provide metrics for probabilities from an information-theory point of view. The idea is to find the shortest description of the data. Assume that the base of the logarithm is 2. The entropy $H(P)=\int_{\mathcal{X}}P(x)\log P(dx)$ describes the optimal average asymptotical code length for a single distribution $P$ according to asymptotic equipartition property~\cite{cover2012elements}. For two distributions $P$ and $Q$ for the same measurable space $\mathcal{X}$, assume that $P$ is the true distribution and $Q$ is the proposed distribution for approximating $P$, the relative entropy, also called Kullback-Leibler divergence, $D(P||Q)=\int_{\mathcal{X}}P(x)\log \frac{P(dx)}{Q(dx)}$ satisfies~\cite{slides1p36:online}
\begin{equation}\label{equ:hd}
    H(P)+D(P||Q) \leq \mathbb{E}L < H(P)+D(P||Q) +1
\end{equation}
where $\mathbb{E}L$ the expected code length for encoding $P$. According to \cref{equ:hd}, if we encode $P$ optimally, $D(P||Q)$ can be interpreted as the extra symbols we need to encode the $Q$. Consequently, minimizing $D(P||Q)$ reduces the redundancy in coding. The maximum likelihood method minimizes $D(P||Q)$ directly. There are many other variants for the divergence metrics. For example, the vanilla GAN proposed in \cite{goodfellow2014generative} is equivalent to minimizing the Jensen-Shannon divergence~\cite{lin1991divergence}
\begin{equation}
    d(P, Q) = \frac{1}{2} D(P||M) + \frac{1}{2} D(Q||M), M = \frac{P+Q}{2}
\end{equation}
 The least square GAN~\cite{mao2017least} is proved to minimize the $\chi^2$ divergence for certain parameters. There is also a class of GANs that minimizes $f$-divergence~\cite{nowozin2016f}, which is also known as Ali-Silvey distances
\begin{equation}
    d(P, Q) = \int_{\mathcal{X}} Q(x)f\left(\frac{P(dx)}{Q(dx)} \right)
\end{equation}
where $f:\mathbb{R}_+\rightarrow \mathbb{R}$, $f(1)=0$ is a convex, lower semi-continuous function.

\textit{Integral probability metrics}~(IPM)~\cite{sriperumbudur2009integral} is a class of distance measures on probabilities that measures the largest discrepancy in expectation over a class of witness functions, which is defined as
\begin{equation}
    d(P, Q) = \sup_{f\in \mathcal{F}}\left| \int fdP - \int fdQ \right|
\end{equation}
where $P$, $Q$ are two probabilities and $\mathcal{F}$ is a class of real-valued bounded measurable functions. Examples of such distances are shown in \cref{tab:ipm}; in the table, $\|f\|_L$ is the Lipschitz semi-norm of a bounded continuous real-valued function $f$ 
\begin{equation}
    \|f\|_L := \sup\left\{ \frac{|f(x)-f(y)|}{|x-y|}: x\neq y \right\}
\end{equation}
\begin{table}[ht]
\centering
\begin{tabular}{@{}ccc@{}}
\toprule
$\mathcal{F}$ & Description  & Note \\ \midrule
$\{f: \|f\|_L \leq 1\}$ & Kantorovich metric & \makecell{The dual representation of\\ the Wasserstein distance} \\
$\{f: \|f\|_L+\|f\|_\infty \leq 1\}$ & Dudley metric &  \\
$\{f: \|f\|_\infty \leq 1\}$ &  \makecell{$L_1$ distance\\ Total variation metric$\times 2$}& \makecell{Equivalent to $\int |P(dx)-Q(dx)|$.\\ The only nontrivial metric that belongs to \\  $f$-divergence and IPM} \\
$\{\mathbf{1}_{(-\infty, t]}: t\in \mathbb{R}^d\}$ & Kolmogorov distance &  \\ 
$\{f: \|f\|_{\mathcal{H}} \leq 1\}$ & Kernel distance & \makecell{$\mathcal{H}$ represents a \\reproducing kernel Hilbert space~(RKHS)} \\\bottomrule
\end{tabular}
\caption{Examples of IPM~\cite{sriperumbudur2012empirical}.}
\label{tab:ipm}
\end{table}

The maximum mean discrepancy belongs to the kernel distance,
\begin{equation}\label{equ:mmd}
    d(P, Q) = \sup_{f\in \mathcal{F}} \left( \mathbb{E}_{x\sim P}f(x) - \mathbb{E}_{y\sim Q}f(y) \right)
\end{equation}
When $\mathcal{F}=C^0(\mathcal{X})$, i.e., the set of all continuous, bounded functions on $\mathcal{X}$~\cite{dudley2018real}, we have $d(P,Q)=0$ if and only if $P=Q$ and thus it can be used as a metric for proximity of $P$ and $Q$. Such functions can be approximated with functions in a universal reproducing kernel Hilbert space~\cite{steinwart2001influence}. The latter metric is used to design MMD GAN~\cite{li2017mmd}.

\textit{Optimal transport}~(OT)~\cite{peyre2019computational} can also be used to derive probability metrics. It basically considers the cost of transforming one probability distribution to another. Let $\mathcal{X}$ and $\mathcal{Y}$ be two measurable space. $P$, $Q$ are two probability distributions on $\mathcal{X}$ and $\mathcal{Y}$ respectively. The Kantorovich problem in optimal transport is
\begin{align}
    \min_{\pi}&\  \int c(x,y) d\pi(x,y) \label{equ:pi} \\
    s.t.&\   P_{\mathcal{X}\#}\pi = P \quad P_{\mathcal{Y}\#}\pi = Q
\end{align}
where $P_{\mathcal{X}\#}\pi$ and $P_{\mathcal{Y}\#}\pi$ denotes the marginals of $\pi$ with respect to the first and second component. If $\mathcal{X}=\mathcal{Y}$, and let $d(x,y)$ be a distance on the space $\mathcal{X}$. If $c(x,y)=d(x,y)^p$, $p\geq 1$ in \cref{equ:pi}, the $p$-Wasserstein distance on $\mathcal{X}$ is defined as
\begin{equation}
    d(P,Q) = \left(\min_{\mu\in \mathcal{L}(P,Q)} \int c(x,y)^p d\pi(x,y)\right)^{1/p}
\end{equation}
where $\mathcal{L}(P,Q)$ is the set of all measures with marginals $P$ and $Q$, and $\rho\geq 0$ is the cost function. For example, if $\mathcal{X}=\mathbb{R}^n$, we can pick $d(x,y)=|x-y|$. Different from many divergences, the Wasserstein distance is a true distance, i.e., $d(x,y)=d(y,x)$ and satisfies the triangle inequality. When $p=1$, the 1-Wasserstein distance has  a supremum form due to Kantorovich-Rubinstein duality
\begin{equation}
    d(P,Q) = \sup_{\|f\|_L\leq 1} \mathbb{E}_{x\sim P}[f(x)] - \mathbb{E}_{x\sim Q}[f(x)]
\end{equation}
where $\|f\|_L\leq 1$ indicates $f$ is 1-Lipschitz. The 1-Wasserstein distance is adopted in Wasserstein GAN~\cite{arjovsky2017wasserstein}.

Good performance in one metric does not necessarily mean the same in another metric. One would pick an appropriate metric for her specific applications. Also the choice of metrics has significant impact on the optimization algorithm. As an illustration, \cref{tab:distance} shows the probability metric for two delta distribution and two Gaussian distributions. If we are comparing discrete distributions and we are using KL-divergence, JS-divergence or the total variation metric, we may get a constant or infinity number, and consequently the gradients will be meaningless. 
\begin{table}[ht]
\begin{tabular}{@{}cc@{}}
\toprule
Metrics & $P=\delta_{\theta_1}$, $Q=\delta_{\theta_2}$  \\ \midrule
KL-divergence & $+\infty$ \\
JS-divergence & $\log 2$\\
Wasserstein-2 distance & $|\theta_1-\theta_2 |$\\
Total variation distance & 1\\
\makecell{IPM \\ $\mathcal{F}=\{x, x^2\}$} & $\max\{|\theta_1, \theta_2 |, |\theta_1^2-\theta_2^2 | \}$ \\ \bottomrule
\end{tabular}
\caption{Different probability metrics for two delta distribution. Here $\theta_1\neq \theta_2$}
\label{tab:distance}
\end{table}

\section*{\refname}

\bibliographystyle{elsarticle-num}
\bibliography{adversarial.bib}

\end{document}